\documentclass[12pt,twoside,final,psamsfonts]{amsart}
\usepackage[psamsfonts]{amssymb}
\usepackage{times,a4wide}
\usepackage[dvips,final]{graphicx,psfrag}
\usepackage[draft=false]{hyperref}
\theoremstyle{plain}
\newtheorem*{theorem*}{Theorem}
\newtheorem{theorem}{Theorem}[section]
\newtheorem*{theoremHL}{Theorem A}

\newtheorem{corollary}[theorem]{Corollary}
\newtheorem{lemma}[theorem]{Lemma}

\theoremstyle{definition}
\newtheorem*{definition}{Definition}
\newtheorem{remark}{Remark}

\newenvironment{proof*}{\vskip 2mm\noindent {}}
\newcommand{\C}{{\mathbb C}}
\newcommand{\N}{{\mathbb N}}
\newcommand{\NN}{{\mathcal N}}
\newcommand{\D}{{\mathbb D}}
\newcommand{\R}{{\mathbb R}}
\newcommand{\T}{{\mathbb T}}

\newcommand{\la}{{\lambda}}
\newcommand{\eps}{{\varepsilon}}

\newcommand{\Hol}{\operatorname{Hol}}
\newcommand{\Har}{\operatorname{Har}}

\newcommand{\De}{\mathcal D_\psi}

\newcommand{\eit}{{e^{i\theta}}}

\renewcommand{\aa}{{A^{-\alpha}}}

\newcommand{\tqnk}{{\tilde Q_{n,k}}}
\newcommand{\liL}{\lambda\in\Lambda}

\renewcommand{\H}{{\mathcal H}}

\author{Xavier Massaneda \& Pascal J. Thomas}

\address{Departament de Matem\`atica Aplicada i An\`alisi,
Universitat  de Barcelona, Gran Via 585, 08007-Bar\-ce\-lo\-na, Spain}

\address{Laboratoire de Math\'ematiques Emile Picard, UMR CNRS 5580,
Universit\'e Paul Sabatier, 118 route de Narbonne, 31062 TOULOUSE
CEDEX, France.}

\thanks{Both authors supported by a XTI program of the Comunitat de
Treball dels Pirineus.
First author also
supported by the DGICYT grant MTM2005-08984-C02-02 and the CIRIT grant
2005-SGR00611.}

\email{xavier.massaneda@ub.edu, pthomas@cict.fr}

\date{\today}


\title[Sampling Sets for the Nevanlinna class] {Sampling Sets for the Nevanlinna class}

\begin{document}

\begin{abstract} We propose a definition of sampling set for the Nevanlinna
and Smirnov classes in the disk and show its equivalence with the notion of determination set for the same classes.
We also show the relationship with determination sets for  related classes of functions
and deduce a characterization of Smirnov sampling sets. For Nevanlinna sampling
we give general conditions (necessary or sufficient), from which we obtain precise geometric
descriptions in several regular cases.
\end{abstract}

\maketitle

\tableofcontents

\section{Introduction}

Let $\Lambda$ be a subset in the unit disk $\D$. In
general  $\Lambda$ is called sampling for a space of
holomorphic functions $X$ when any function $f\in X$ is determined by its
restriction $f|\Lambda$, with control of
norms. For Banach spaces $X$ it is usually clear what that control of norms
means, but for the spaces we have in mind the situation is not so obvious. Consider the
\emph{Nevanlinna class}
\[
\NN=\bigl\{f\in \Hol(\D):\lim_{r \to 1}\frac{1}{2\pi}\int_{0}^{2\pi}
     \log^+|f(re^{i\theta})|\;d\theta<\infty\bigr\}\ ,
\]
which is not a Banach space, but  enjoys the structure of
complete metric space with the distance $d(f,g)=N(f-g)$ induced by
\[
N(f)=\lim_{r \to 1}\frac{1}{2\pi}\int_{0}^{2\pi}
     \log(1+|f(re^{i\theta})|)\;d\theta\ .
\]
The subharmonicity of $\log (1+|f|)$ yields the pointwise estimate
\[
(1-|z|)\log (1+|f(z)|)\leq 2N(f),
\]
which shows that convergence in the distance $d$ implies 
uniform convergence on compact sets \cite[Proposition 1.1]{SS}.

The Nevanlinna class $\NN$ coincides with the set of functions $f\in\Hol(\D)$ such
that $\log(1+ |f|)$ admits a harmonic majorant (see \cite[p.69]{Gar} or \eqref{deco} later on).
The value $N(f)$ can then be rewritten in terms of an extremal problem for 
harmonic majorants. Let $\Har_+(\D)$ denote the space of non-negative harmonic functions in the
disk; then
\[
N(f)=\inf\{h(0)\, :\, \textrm{$h\in\Har_+(\D)$ with $\log(1+|f|)\leq h$}\}\ .
\]
This expression makes sense for any $f$ measurable in $\D$, in particular for a restriction $f|\Lambda$, and
suggests the following definition.

\begin{definition}
A set $\Lambda$ is \emph{sampling for} $\NN$ if there exists $C>0$ such that
\[
N(f)\leq N(f|\Lambda)+ C\qquad\forall f\in \NN\ .
\]
\end{definition}

In Section \ref{S1} we study
first the relationship between Nevanlinna sampling sets and determination sets for the same class and for 
the space $\H^\infty$ of bounded holomorphic functions. We prove that sampling and 
determination sets for $\NN$ are the same. Also, from the characterization of 
$\H^\infty$ determination sets given by Brown, Shields and Zeller \cite{BrShZe} we
deduce a complete
description of sampling sets for the 
\emph{Smirnov class}
\[
     \NN^+=\bigl\{f\in \NN:\lim_{r \to 1}\frac{1}{2\pi}\int_{0}^{2\pi}
     \log^+|f(re^{i\theta})|\;d\theta=\frac{1}{2\pi}\int_{0}^{2\pi}
     \log^+|f^*(e^{i\theta})|\;d\theta<\infty\bigr\}.
\]
Here $f^*(e^{i\theta})$ denotes the  non-tangential limit of $f$ at the
boundary point $e^{i\theta}$.

Next we study the relationship between sampling sets for $\NN$ and determination 
sets for the class
$\Har_\pm(\D)$ of harmonic functions which can be written as the difference 
of two positive harmonic functions. This is relevant because  the 
Riesz-Smirnov factorization implies that  for any $f\in\NN$ there exist a Blaschke
product $B$ and $h\in\Har_{\pm}(\D)$ such that $\log|f|=\log|B|+h$. From the
characterization
of determination sets for $\Har_{\pm}(\D)$ given by Hayman and Lyons \cite{HaLy} we deduce a necessary
geometric condition for sampling in $\NN$. Some examples show that this condition is far from
being sufficient.

In Section \ref{S2} we give  general  conditions for Nevanlinna
sampling (Theorem~\ref{suficient}), which in Section \ref{S3} are used
to obtain a precise geometric description for three different
types of regular sampling sets: fine nets of points, regular sequences on cercles tending to
the unit cercle, and uniformly dense unions of hyperbolic disks, as
considered by  Ortega-Cerd\`a and Seip in \cite{OrSe}. 
%It is interesting to
%note that the condition for  sampling in this situation is the same as the
%condition in \cite[Theorem1]{OrSe}.

A final remark about notation. The expression $A\preceq B$ means that there
exists a constant $C>0$, independent of whatever arguments are involved, such
that $A\leq C B$. If both $A\preceq B$ and $B\preceq A$ then we write $A\simeq
B$.

\section{Determination sets and necessary conditions}\label{S1}

In this section we describe the relationship between our definition of sampling 
and other related notions studied previously. 
%From these we deduce a necessary 
%condition (Theorem \ref{HayLyons}).

\subsection{Sampling and determination sets}
We begin with an easy observation: in the definition of sampling given in the introduction
$N(f)$ can be equivalently replaced by 
\[
N_+(f)=\lim_{r\to 1}\frac 1{2\pi}\int_0^{2\pi}\log^+|f(r e^{i\theta})|d\theta=
\inf\{h(0)\, :\, \textrm{$h\in\Har_+(\D)$ with $\log^+|f|\leq h$}\}\ .
\]
This is an immediate consequence of the inequalities
\[
\log^+|f|\leq \log(1+|f|)\leq\log 2+\log^+|f|\ .
\]

Also, the constant $C$ can be assumed to be 0, as the following Lemma shows.

\begin{lemma}\label{constant} A set $\Lambda $ is sampling for $\NN$ if and only if
$N_+(f)=N_+(f|\Lambda)$ for all $f\in \NN$.

\end{lemma}

\begin{proof}
Of course, we only need to see that the equality is necessary. Since $\NN$ is an algebra, the
sampling inequality $N_+(f)\leq N_+(f|\Lambda)+ C$ yields automatically
\[
N_+(f^n)\leq N_+(f^n|\Lambda)+C\qquad\forall f\in \NN\quad\forall n\in\N\ .
\]
By definition 
$N_+(f^n)=n N_+(f)$, so
\[
N_+(f)\leq N_+(f|\Lambda)+\frac Cn\qquad\forall f\in \NN\quad\forall n\in\N\ ,
\]
and the result follows letting $n$ tend to $\infty$.
\end{proof}

Let us consider also two related notions for a set $\Lambda$ in
$\D$. A priori, one seems weaker and the other  stronger than the sampling
property. 

%Let $\H^\infty$ denote the space of bounded holomorphic
%functions in the unit disk.

\begin{definition}
A set $\Lambda$ is  a \emph{determination set for} $\NN$ if
$\NN\cap L^\infty(\Lambda)\subset \H^\infty$, i.e. if
any $f\in \NN$ with $\sup_\Lambda |f|<\infty$ must be
bounded on the whole unit disk.

A set $\Lambda$ is  \emph{strongly sampling for} $\NN$ if
whenever $f\in \NN$ and $h\in\Har_+(\D)$ are such that $\log^+|f(\lambda)|\leq
h(\lambda)$ for all $\liL$, then necessarily $\log^+|f|\leq h$.

When  $\Lambda$ is  strongly sampling the distance $d(f,g)$ between two functions
$f,g\in\NN$ coincides with the distance between their restrictions $f|\Lambda$ and
$g|\Lambda$.
\end{definition}

\begin{remark} \label{BSZ} 
A set $\Lambda$ is a \emph{determination set for} $\H^\infty$ 
when $\|f\|_\infty=\sup_\Lambda |f|$ for all $f\in \H^\infty$. It is easy to see
that determination sets for $\NN$ are also determination sets for $\H^\infty$, 
which therefore satisfy $\|f\|_\infty=\sup_\Lambda |f|$ for all $f\in\NN$.

Indeed, assume that there exists $f\in \H^\infty$ such that $\|f\|_\infty=1$
and $ \sup_\Lambda |f|=s<1$. Take $\{z_k\}_k\subset\D$ such that
$\lim_k |f(z_k)|=1$ and consider any accumulation point $\zeta\in\T$ of
$\{f(z_k)\}_k$. Then the function $g=1/(\zeta-f)$ belongs to $\NN$, is not bounded,
but $\sup_\Lambda|g|\leq 1/(1-s)$.

Brown, Shields and Zeller showed that
$\Lambda$ is a determination set for $\H^\infty$  if and only if
the set $NT(\Lambda)$ consisting of the $\zeta\in\T$ which
are a non-tangential limit of points in $\Lambda$
has full measure, i.e. $|NT(\Lambda)|=2\pi$  \cite{BrShZe}. It was
shown in \cite{Th} that the same condition also characterizes
sampling sets for the Hardy spaces $\H^p$ ($0<p<\infty$), if appropriately defined. 
This condition
is therefore necessary for $\Lambda$ to be a determination set for $\NN$.
\end{remark}

Our first result shows that the previous notions are all
equivalent.

\begin{theorem}\label{equivalencia}
Let $\Lambda$ be a subset of $\D$.
The following properties are equivalent:

\begin{itemize}
\item[(a)] $\Lambda$ is a sampling set for $\NN$.
\item[(b)] $\Lambda$ is a determination set for $\NN$.
\item[(c)] $\Lambda$ is a strongly sampling set for $\NN$.
\end{itemize}
\end{theorem}

It is clear from (c) that the sampling property is invariant
under automorphisms of the disk: if  $\Lambda$ is sampling for
$\NN$ and $\phi(z)=e^{i\theta}\frac {z-a}{1-\bar a z}$, $a\in\D$, is an
automorphism of $\D$, then $\phi(\Lambda)$ is also sampling for $\NN$.

Before the proof we need to recall some well-known facts about the structure
of the Nevanlinna  class (general refe\-rences are e.g. \cite{Gar},
\cite{Nik02} or \cite{RR}). 

For a set $Z \subset \D$ with multiplicities, the  \emph{Blaschke
product} with zeros on $Z$ is
$$
B_Z (z) := \prod_{a\in Z} \frac{\bar a}{|a|}\frac{a-z}{1-z\bar a},
$$
where the points are repeated according to multiplicities.  This is convergent,
not identically equal to $0$, if and only if $\sum_{a\in Z} (1-|a|) < \infty$.
When this is the case, we say that $Z$ is a \emph{Blaschke sequence}, or verifies the
Blaschke condition.

A function $f$ is called \emph{outer} if it can be written in the form
\[
f(z)=C \exp \left\{\int_0^{2\pi} \frac{e^{i\theta}+z}{e^{i\theta}-z}
\log v(e^{i\theta}) \frac{d\theta}{2\pi}\right\},
\]
where $|C|=1$, $v>0$  a.e.\ on $\T$ and $\log v \in L^1(\T)$. Such a
function is
the quotient $f=f_1/f_2$ of two bounded outer functions  $f_1,f_2\in \H^\infty$
with  $\|f_i\|_{\infty}\le 1$, $i=1,2$. In particular, the weight $v$
is given by
the boundary values of $|f_1/f_2|$. Setting $w=\log v$,
%and writing $P[w]:=$P[wd\sigma]$,
we have
\[
\log |f(z)|=P[w](z) :=\int_0^{2\pi} P_z(e^{i\theta}) 
w(e^{i\theta})\frac{d\theta}{2\pi}\ ,
\]
where
$$
P_z(e^{i\theta}) := \frac{1-|z|^2}{|e^{i\theta} - z|^2}
$$
denotes the \emph{Poisson kernel} at $z \in \D$.

In general, for any finite measure $\mu$ on $\T$, the \emph{Poisson integral}
of $\mu$ is the harmonic function given by
$$
P[\mu](z) := \int_{0}^{2\pi} P_z(e^{i\theta}) \, d\mu(\theta).
$$

Another important family in this context are  \emph{inner} functions: $I\in
\H^\infty$ such that $|I|=1$ almost everywhere on $\T$. Any inner function $I$ can be
factorized into a Blaschke product $B$ carrying the zeros
of $I$, and a singular inner function
$S$ defined by
\begin{equation*}
S(z)=\exp\left\{-\int_0^{2\pi} \frac{e^{i\theta}+z}{e^{i\theta}-z}\,
d\mu(e^{i\theta})\right\},
\end{equation*}
for some positive Borel measure $\mu$ singular with respect to Lebesgue
measure.

According to the Riesz-Smirnov factorization, any
function $f\in \NN$ is represented as
\begin{equation}\label{fact}
f=\alpha \frac{B S_1 f_1}{S_2 f_2},
\end{equation}
with $f_i$ outer, $\|f_i\|_{\infty}\le 1$, $S_i$ singular inner,
$B$ a Blaschke product and $|\alpha|=1$.

\begin{remark}\label{factor}
Let $\Har_{\pm}(\D)$ denote the set of harmonic functions
$h$ that can be written $h=h_1-h_2$, with $h_1,h_2\in\Har_+(\D)$. 
The factorization above shows that for $f\in\NN$ there exist always
$h\in\Har_{\pm}(\D)$  and a Blaschke product $B$ such that
\begin{equation}\label{deco}
\log|f|=\log|B|+h\ ,
\end{equation}
and reciprocally, for any $h\in\Har_{\pm}(\D)$  and any Blaschke product $B$
there exists $f\in\NN$ satisfying \eqref{deco}.
\end{remark}

\begin{proof}[Proof of Theorem~\ref{equivalencia}] (c)$\Rightarrow$(a) is
immediate from the definition.

(a)$\Rightarrow$(b). Let $f\in \NN$ with $s_\Lambda=:\sup_\Lambda|f|<\infty$ and
consider $g=f/s_\Lambda\in \NN$. Since $\Lambda$ is sampling and $\log^+|g(\lambda)|=0$
for all $\liL$ we have, according to Lemma~\ref{constant},  $N_+(g)=N_+(g|\Lambda)=0$. 
%If $h\in\Har_+(\D)$ is such that $\log^+|g|\leq h$ we have
%\[
%\int_0^{2\pi}\log^+|g(re^{i\theta})|\, d\theta\leq
%\int_0^{2\pi}h(re^{i\theta})\, d\theta=h(0)\qquad r\leq 1 .
%\]
Thus 
%$N_+(g)=0$ implies that
$\int_0^{2\pi}\log^+|g(re^{i\theta})|\, d\theta=0$ for all $r<1$, hence 
%by the maximum principle
$\|g\|_\infty\leq 1$.

(b)$\Rightarrow$(c). Let $f\in \NN$ and $h\in\Har_+(\D)$ be such that
$\log|f(\lambda)|\leq h(\lambda)$ for all $\liL$. By Remark~\ref{factor},
there exists a function $g\in \NN$ such that $\log|g|=\log|f|-h$. 
We have then $\log |g(\lambda)|\leq 0$ for all $\liL$, and as 
pointed out in Remark~\ref{BSZ}, this implies
 $\|g\|_\infty \leq 1$,  i.e.  $\log|g|=\log|f|-h\leq 0$.
\end{proof}

\subsection{Sampling in the Smirnov class}
All the definitions and proofs above can be similarly given for the
Smirnov class $\NN^+$ defined in the introduction. The Smirnov class
consists of those $f\in \NN$ for which the harmonic
majorant of $\log^+|f|$ is
quasi-bounded (the Poisson integral of some $w\in L^1(\T)$). 
Equivalently, it consists of those $f\in \NN$ with no singular factor 
$S_2$ in the factorization \eqref{fact}.

The geometric description of sampling sequences for $\NN^+$ is a 
straightforward consequence of the results in \cite{BrShZe} and
Remark~\ref{BSZ}. Recall that $NT(\Lambda)$ denotes the 
non-tangential accumulation set of $\Lambda$ in $\mathbb T$.

\begin{theorem}\label{smirnov}
Let $\Lambda$ be a subset of $\D$.
The following  properties are equivalent:

\begin{itemize}
\item[(a)] $\Lambda$ is a sampling set for $\NN^+$.
\item[(b)] $\Lambda$ is a determination set for $\NN^+$.
\item[(c)] $\Lambda$ is a strongly sampling set for $\NN^+$.
\item[(d)] $|NT(\Lambda)|=2\pi$.
\end{itemize}
\end{theorem}

\begin{proof} The equivalence between (a), (b) and (c) is seen as in
Theorem~\ref{equivalencia}.

The necessity of (d) is pointed out in Remark~\ref{BSZ}. The
sufficency is immediate: for almost every
$\theta\in [0,2\pi)$ there exists a sequence $\{\lambda_k\}_k\subset\Lambda$
tending non-tangentially to $e^{i\theta}$, and therefore
$f^*(e^{i\theta})=\lim_{k\to\infty} f(\lambda_k)$ 
\cite[Theorem 5.3]{Gar}. Then, if $f\in
\NN^+$ and $h\in\Har_+(\D)$ are such that $\log^+|f(\lambda)|\leq h(\lambda)$
for all $\liL$ we have $\log^+|f^*(e^{i\theta})|\leq h(e^{i\theta})$ a.e
$\theta\in\T$. This yields $N_+(f)\leq N_+(f|\Lambda)$.
\end{proof}

\subsection{Determination sets for harmonic functions and a necessary 
condition for Nevanlinna sampling}\label{nec}

From previous results on determination sets for harmonic functions and the
equivalences of Theorem~\ref{equivalencia} we deduce a first necessary
condition for Nevanlinna sampling (Corol\-la\-ry~\ref{HLnec}). 
This can be obtained directly, as shown in the Appendix.

Given $z,w\in\D$ let 
\[
\rho (z,w):= \left| \frac{z-w}{1-z\bar w}\right|
\]
stand for the
its \emph{pseudohyperbolic distance}. For $r\in (0,1)$ and $z\in\D$ let
$D (z,r)=\{w\in\D : \rho(z,w)<r\}$.

A sequence $\Lambda=\{\lambda_k\}_k$ is called \emph{separated} when
\[
\inf_{j\neq k} \rho(\lambda_j,\lambda_k)>0\ .
\]

For any set $\Lambda \subset \D$  and $\delta \in (0,1)$, consider
the pseudohyperbolic dilation
$$
\Lambda^\delta = \bigcup_{\lambda \in \Lambda} D(\lambda,\delta),
$$
and given $\zeta \in \T$ denote
$$
I(\Lambda,\zeta, \delta) := \int_{\Lambda^\delta} 
\frac1{|\zeta - z|^2} dm(z)\ ,
$$
where $dm$  stands for the usual area measure.

We note that for any fixed $\zeta \in \T$, the values 
$I(\Lambda,\zeta, \delta) $
are finite simultaneously for all values of $\delta \in (0,1)$, and that
this is equivalent to the fact that for any maximal separated subsequence
$\Lambda'\subset\Lambda$, we have
$$
\sum_{\lambda\in  \Lambda'} \bigl( \frac{1-|\lambda|^2}{|\zeta - \lambda|} 
\bigr)^2
= \sum_{\lambda\in \Lambda'}  (1-|\lambda|^2) P_\lambda (\zeta)  < \infty.
$$

We recall the following characterization of determination sets for the class 
$\Har_\pm(\D)$ given by Hayman and Lyons \cite{HaLy}. This is elaborated upon in
\cite{Ga}. 

\begin{theoremHL}\label{HayLyons}
Let $\Lambda \subset \D$.
The following properties are equivalent.
\begin{itemize}
\item[(a)] $\sup_{\Lambda} h  = \sup_{ \D} h\ $ 
for all $h \in \Har_{\pm}(\D)$.
\item[(b)]
There exists $ \delta \in (0,1)$ such that
$I(\Lambda,\zeta, \delta) = \infty\ $ for all $\zeta \in \T$.
\end{itemize}
\end{theoremHL}

We shall call the sets satisfying these condition \emph{Hayman-Lyons} sets. 

Condition (b) is more restrictive than Brown, Shields and Zeller's condition $|NT(\Lambda)|=2\pi$.
Actually, $|NT(\Lambda)|=2\pi$ is equivalent to $I(\Lambda,\zeta, \delta) = \infty\ $ a.e. $\zeta \in \T$
\cite[Corollary 2]{Ga}. On the other hand, it is clear that if $NT(\Lambda)=\mathbb T$
then (b) is satisfied, since the Poisson kernel $P(z,\zeta)$ is
bounded below in any Stolz angle with vertex at $\zeta$.

\begin{corollary}\label{HLnec}
\label{neccond}
A Nevanlinna sampling set is a Hayman-Lyons set.
\end{corollary}

\begin{proof}
By Theorem~\ref{equivalencia}, $\Lambda$ is a determination set for $\NN$, and therefore
$\sup_\Lambda\log |f|=\sup_{\D}\log |f|$ for all $f\in\NN$. By Remark~\ref{factor}, 
this implies (a) in Theorem A.
\end{proof}

Notice also that when $\Lambda$ is a determination set for $\NN$
and $f,g\in\NN$ are such that $|f(\lambda)|\leq |g(\lambda)|$ for
all $\lambda\in\Lambda$, then $|B f|\leq |g|$, where $B$ indicates
the Blaschke product associated to the zeros of $g$. To see this
factorize $g=Bg_0$, with $g_0$ non-vanishing. Then
$|f(\lambda)|/|g_0(\lambda)|\leq |B(\lambda)|\leq 1 $ and by
hypothesis $|f|\leq |g_0|$, as desired.

The Hayman-Lyons condition is not sufficient for sampling in
$\NN$, as shown in the following example.

{\bf{Example 1.}}
Take a dyadic partition of the
disk: for any $(n,k)$ in the set of indices
$\mathcal I=\{(n,k) : n\in\mathbb N,\; 0 \le k \leq  2^n-1\}$
consider the interval
\begin{equation}
\label{dyadarcs}
I_{n,k} := \{ \eit : \theta \in [{2\pi}k 2^{-n},{2\pi}(k+1)
2^{-n}) \},  
\end{equation}
and the associated Whitney partition in ``dyadic squares":
\begin{equation}
\label{dyadsquares}
Q_{n,k} := \{ r \eit : \eit \in I_{n,k} , 1- 2^{-n} \le r < 1-
2^{-n-1} \}.
\end{equation}
Observe that the pseudohyperbolic diameter of each Whitney square
$Q_{n,k}$ is bounded between two absolute constants.

Let $\Lambda$ be the sequence consisting of the centers 
$c_{n,k}$ of $ Q_{n,k}$. An immediate computation shows that 
for every $\zeta\in\T$
\[
\sum_{n=1}^\infty\sum_{k=0}^{2^n-1}(1-|c_{n,k}|^2) P_{c_{n,k}} 
(\zeta) =
\sum_{n=1}^\infty\sum_{k=0}^{2^n-1} \bigl( 
\frac{1-|c_{n,k}|^2}{|\zeta - c_{n,k}|} \bigr)^2
\simeq \sum_{n=1}^\infty\sum_{k=0}^{2^n-1} 
\bigl(\frac{2^{-n}}{2^{-n}+k 2^{-n}} \bigr)^2
=  \infty,
\]
and therefore $\Lambda$ is a Hayman-Lyons set. 

In order to see that $\Lambda$ is not a determination set for $\NN$ 
fix $\zeta=1\in\T$ and consider a \emph{horocycle}  $\{z: P_z(1)= 
c\}$, the boundary of the euclidian disk $ B(\frac c{1+c}, \frac 1{1+c})$
 (see Figure 1). Then 
$Z=\Lambda\cap B(\frac c{1+c}, \frac 1{1+c})$ is a Blaschke sequence:
\[
\sum_{a\in Z} 1-|a|\simeq \sum_{n=1}^\infty\sum_{0\leq
k<\sqrt{2^n}}(1-|c_{n,k}|)\simeq \sum_{n=1}^\infty 2^{-n/2}<\infty \ .
\]

\vglue 0.5 truecm
\begin{center}
%\psfrag{tag0}{$0$}
\psfrag{tag1}{\small{$Q_{n,k}$}}
\psfrag{tag2}{\small{$c_{n,k}$}}
\psfrag{tag3}{\footnotesize{$P_z(1)=c $}}
\psfrag{tag4}{\small{$\zeta $}}
%\psfrag{...}{$\cdots$}
\includegraphics[totalheight=7cm]{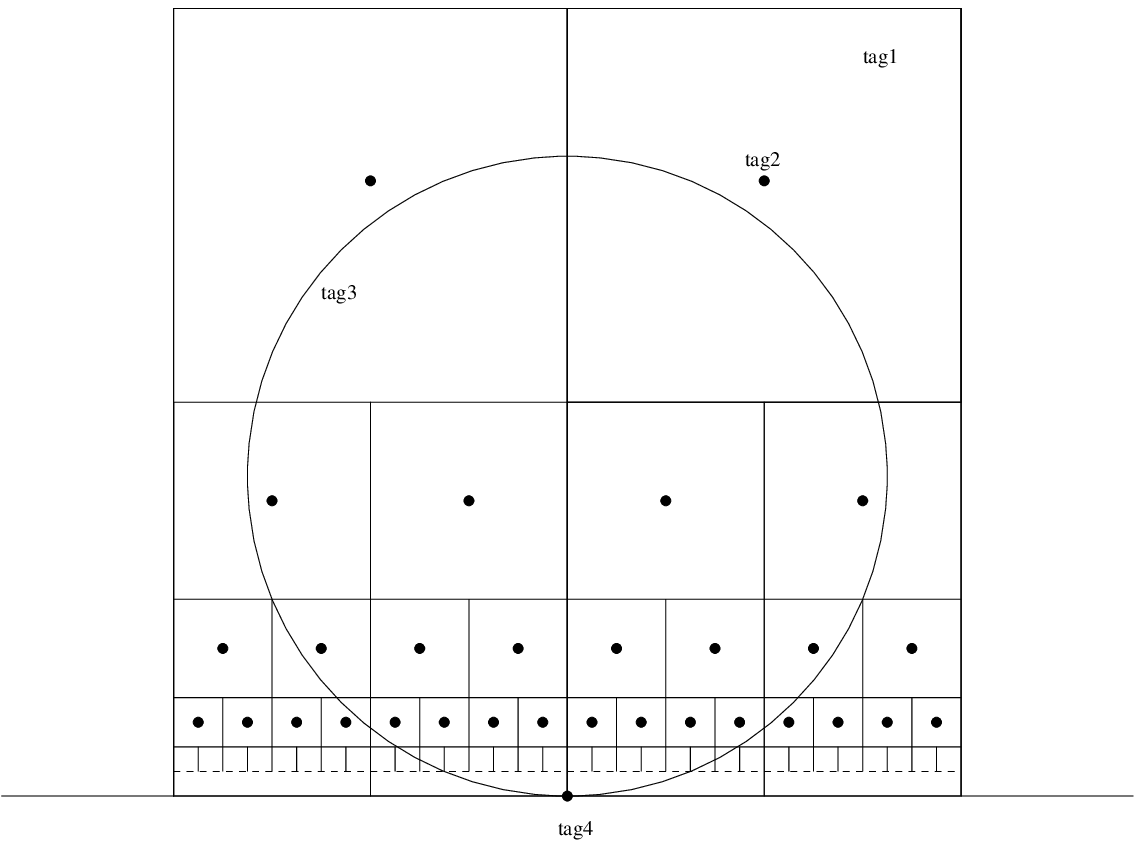}

{\footnotesize{Figure 1. Upper half-plane representation of $\Lambda$ and the horocycle $P_z(1)=c$.}}
\end{center}
\vglue 0.5 truecm

Therefore, there exists $f\in\NN$ such that
\[
\log |f(z)|= \log |B_Z(z)|+P_z(1)\qquad z\in\D\ .
\]

Clearly $\log |f(\lambda)|\leq  c $ for all $\lambda\in\Lambda$.
The fact that $f$ is not bounded rests on the following more general lemma.

\begin{lemma}\cite[Lemma 2.3]{SS}\label{blaschke}
\label{blasgrowth}
For any Blaschke product $B$ and any $\zeta \in \T$,
$$
\limsup_{r\to 1} (1-r) \log |B(r\zeta)| = 0.
$$
\end{lemma}

This implies in particular that 
$\limsup_{r\to 1} (1-r) \log|f(r)| =\limsup_{r\to 1} (1-r) P_r(1)=2$, so $f$
cannot be bounded on $(0,1)$.

\section{General conditions}\label{S2}

In order to see what extra conditions are required on a Hayman-Lyons set 
$\Lambda$ to be a determination set for 
$\NN$, assume that $f\in \NN$ is
such that $\sup_\Lambda|f|\leq 1$. According to Remark~\ref{factor},
there exist a Blaschke product $B$ (with
zero-sequence $Z$) and $F\in\Har_{\pm}(\D)$ such that
$\log |f| = \log|B|+F$. 

It will enough to see that $F$ is quasi-bounded, that is, that $f$ belongs to the 
Smirnov class. This is so because the Hayman-Lyons condition implies $|NT(\Lambda)|=2\pi$, and
we deduce then from Theorem~\ref{smirnov} that 
$\sup_\Lambda |f|=\sup_{\D}|f|$ (see Remark~\ref{BSZ}).

The hypothesis is, in these terms,
\[
F(\lambda)\leq\log\frac 1{|B(\lambda)|}\qquad \liL\ ,
\]
and we would like to impose certain conditions on $\Lambda$ so that this estimate implies that $F$
has a quasi-bounded harmonic majorant. 

A first observation is that 
the zeros of $B$ far from a given $\lambda$ are no obstruction to such majorization.
The following is a restatement of \cite[Proposition 4.1, pp. 13--14]{HMNT},
and of part of its proof.

\begin{lemma} \label{Artur}
Let $B$ be the Blaschke product associated to a Blaschke sequence $Z$.  
For any $\delta \in (0,1)$, there
exists a positive quasi-bounded  harmonic function $H_B =P[w]$,
$w \in L^1(\T)$, such that
$-\log |B(z)| \leq H_B (z) - \sum_{a\in Z\cap D(z,\delta)} \log\rho(z,a)$,  for any $z \in
\D$.
Furthermore 
$$
w(\zeta) = c_0 \sum_{a \in Z} \chi_{I_a} (\zeta),
$$
where $c_0$ is an appropriate positive
constant, and 
$I_a = \{ \zeta \in \T : |\zeta - \frac{a}{|a|}| \le 1-|a| \}$ 
is the ``Privalov shadow" of $a$ on $\T$.
\end{lemma}

Given $\delta\in (0,1)$,
there exists thus $H_{B}$ harmonic, quasi-bounded and positive such that,
\begin{equation}\label{basica}
F(\lambda)\leq H_{B}(\lambda)+
\sum_{a\in Z\cap D(\lambda,\delta)}\log\frac 1{\rho(\lambda,a)}\qquad \liL\ .
\end{equation}
It is clear that we cannot expect to bound the local sum 
in the right hand side of this inequality
by a quasi-bounded 
harmonic function for all $\liL$ (it could happen, for example, 
that $Z\cap\Lambda\neq\emptyset$). Rather, we would like to find conditions on $\Lambda$ 
that ensure such a bound for a subset $\Lambda'\subset\Lambda$ big enough so that
the estimate of $F$ by a quasi-bounded harmonic function
on $\Lambda'$ implies the same estimate everywhere (in the spirit of the
Hayman-Lyons condition for functions in $\Har_{\pm}(\D)$).

For that purpose we need a measure of the ``vulnerability" of $\Lambda$ to the presence
of zeros of a Blaschke product.

Consider the dyadic squares $\{Q_{n,k}\}_{(n,k)\in\mathcal I}$ defined in \eqref{dyadsquares} and
denote by $\tilde Q_{n,k}$ the union of all $Q_{m,j}$ such that
$\overline{Q_{n,k}} \cap \overline{Q_{m,j}} \neq \emptyset$.
There exists  $\delta_0 >0$ such that
%\begin{equation*} \label{fatsquares}
$\overline{Q_{n,k}}^{\delta_0} \subset \tqnk$.
%\end{equation*}

To measure the vulnerabilty of $\Lambda$ at each $Q_{n,k}$,
for $N\in \N$ and $\delta \in (0,1)$ consider 
$$
w_{n,k}(\Lambda, N) =
\sup \left\{
\left(
\inf_{\lambda\in \Lambda\cap\overline{Q_{n,k}}} \sum_{j=1}^N \log \frac1{\rho(\lambda,a_j)}
\right), a_1, \dots , a_N \in \overline{Q_{n,k}}^\delta
\right\}.
$$
We take an empty sum to be $0$, so that $w_{n,k }(\Lambda,0)=0$ for any set $\Lambda$.
Clearly, $w_{n,k}(\Lambda,N)$ is an increasing function of $N$, and there exists 
$C(\delta)>0$ such that $w_{n,k}(\Lambda,N) \ge N C(\delta)$.

Given a Blaschke sequence $Z$ let $N_{n,k}=\#(Z\cap\tqnk)$. The Blaschke condition is
thus equivalent to $\sum_n2^{-n}\sum_{k=0}^{2^n-1} N_{n,k}<\infty$. Any sequence of integers 
$\{N_{n,k}\}_{(n,k)\in\mathcal I}$ satisfying this condition will be called a
\emph{Blaschke distribution}.

%Consider the normalized Poisson kernels $\tilde P_z(\zeta)=(1-|z|^2)P_z(\zeta)$, 
%and let $c_{n,k}$ denote the center of $Q_{n,k}$.

\begin{theorem}\label{suficient}
Let $\Lambda\subset \mathbb D$. 
Each of the following properties implies the next one.

\begin{itemize}
\item[(a)] For any Blaschke distribution $\{N_{n,k}\}_{(n,k)\in\mathcal I}$, there 
exists $\mathcal Q \subset \mathcal I$ such that
\begin{equation}\label{caixes} 
\Lambda \cap Q_{n,k} \neq \emptyset\quad \textrm{for any $(n,k) \in \mathcal Q$},
\end{equation} 
\begin{equation}\label{HLzeta}
\sum_{(n,k) \in \mathcal Q} 2^{-n}  P_{c_{n,k}}(\zeta) = \infty
\end{equation}
for all $\zeta \in \T$, and 
\begin{equation} \label{invul}
\sum_{(n,k) \in \mathcal Q} 2^{-n} w_{n,k}(\Lambda, N_{n,k}) < \infty .
\end{equation}

\item[(b)]
For any Blaschke distribution $\{N_{n,k}\}_{(n,k)\in\mathcal I}$ and any positive finite 
measure $\nu$ on $\T$, singular with respect to the Lebesgue 
measure, there exists $\mathcal Q \subset 
\mathcal I$ satisfying \eqref{caixes}, \eqref{invul} and
condition \eqref{HLzeta}  almost everywhere with respect to $\nu$. 

\item[(c)]
$\Lambda$ is a determination set for the Nevanlinna class.

\item[(d)]
For any Blaschke distribution $\{N_{n,k}\}_{(n,k)\in\mathcal I}$ and any 
positive finite measure $\nu$ on $\T$, singular with respect to the Lebesgue 
measure, there exists $\mathcal Q \subset  \mathcal I$ satisfying \eqref{caixes}, 
\eqref{invul}  and
\begin{equation}\label{HLint}
\int_{\T} \sum_{(n,k) \in \mathcal Q}  2^{-n} 
  P_{c_{n,k}}(\zeta) d\nu(\zeta) 
= \infty .
\end{equation}

\item[(e)]
For any Blaschke distribution $\{N_{n,k}\}_{(n,k)\in\mathcal I}$ and any 
$\zeta \in \T$,  there exists 
$\mathcal Q \subset \mathcal I$ satisfying \eqref{caixes}, 
 \eqref{HLzeta} and \eqref{invul}.
\end{itemize}

\end{theorem}

As pointed out before the statement of Theorem A, condition \eqref{HLzeta} is precisely
the Hayman-Lyons condition for the set $\Lambda\cap\{Q_{n,k}\}_{(n,k)\in\mathcal Q}$.

We will see in the next Section how this somewhat cumbersome conditions can be used to provide
precise geometric conditions, at least when then set $\Lambda$ has some regularity.

Open question: are there examples of sets 
$\Lambda$ to show that the first and last of 
those properties are actually different ?

\begin{proof}
$\textrm{(a)}\Rightarrow \textrm{(b)}$. Obvious.

$\textrm{(b)}\Rightarrow \textrm{(c)}$. Start, as in the general scheme, with $f\in\NN$ 
such that $\sup_\Lambda|f|\leq 1$ and consider the
decomposition $\log |f|=\log |B| + F$.
By the Riesz-Smirnov factorization \eqref{fact}, the function $F$ can be written as
\[
F=h_1-h_2+H_1-H_2 ,
\]
where $h_i,H_i\in\Har_+(\D)$, $H_i$ are quasi-bounded and $h_i=P[\nu_i]$, 
with $\nu_i$ positive finite measure in $\T$, singular with respect to the
Lebsesgue measure.

In order to see that $f\in\NN^+$ it will be enough to prove that $h_1$ has a 
quasi-bounded majorant. To do that we use (b) with the singular measure $\nu_1$
and the Blaschke distribution determined by $B$.
Let $\mathcal Q\subset\mathcal I$ be the set of indices for which (b) holds.

Let $g\in\NN$ be such that $\log|g|=\log|f|-H_1=\log|B|+h_1-h_2-H_2$ 
(explicitly $g=fe^{-(H_1+i\tilde H_1)}$,
where $\tilde H_1$ denotes the harmonic conjugate of $H_1$). Then, the corresponding 
estimate \eqref{basica} for $g$ becomes
\[
h_1(\lambda)-h_2(\lambda)-H_2(\lambda)\leq H_{B}(\lambda)+
\sum_{a\in Z\cap D(\lambda,\delta)}\log\frac 1{\rho(\lambda,a)}\qquad \liL\ ,
\]
for some positive quasi-bounded harmonic function $H_B$.

For each $(n,k)\in\mathcal Q$  there is a particular 
$j=j(n,k)$ such that  $\la_{j(n,k)} \in \Lambda
\cap \overline{Q_{n,k}}$ and 
\[
\sum_{a \in Z \cap Q_{n,k}^\delta} \log
\frac1{\rho(\la_{j(n,k)}, a)}=
\inf_{\la_j \in \Lambda\cap \overline{Q_{n,k}}}
\sum_{a \in Z \cap Q_{n,k}^\delta} \log
\frac1{\rho(\la_j,a)}
\le w_{n,k}(\Lambda,N_{n,k}).
\]
Let $I_{n,k}$ be the dyadic arcs defined in \eqref{dyadarcs}
and let $C>0$. According to \eqref{invul} the function
$H_L=P [W]$, with 
\[
W:=C\sum_{(n,k)\in\mathcal Q} 
w_{n,k}(\Lambda,N_{n,k}) \chi_{I_{n,k}}\ ,
\] 
is a positive quasi-bounded harmonic function. 
The usual estimate of the 
Poisson kernel by a ``square" kernel (or a direct computation)
shows that
$P[\chi_{I_z}](z) \ge c >0$, with $c$ independent of $z$, hence for
$C$ well chosen,
\[
w_{n,k}(\Lambda,N_{n,k})\leq H_L(\la_{j(n,k)}),
\qquad (n,k)\in\mathcal Q\ .
\] 
Then, defining
$\Lambda'=\{\la_{j(n,k)}\}_{(n,k)\in\mathcal Q}$  we have
$h_1- h_2-H_2\leq H_B+H_L$ on $\Lambda'$.
Condition \eqref{HLzeta} being satisfied $\nu_1$-a.e. and  
\cite[Theorem 2]{Ga} show then that
\[
1\leq \inf_{\Lambda'} \frac{h_2+H_2+H_B+H_L}{h_1}=
\inf_{\D} \frac{h_2+H_2+H_B+H_L}{h_1}\ .
\]
Thus $h_1- h_2-H_2\leq H_B+H_L$ everywhere, and
$\log|f| \leq \log|B|+H_1+H_B+H_L$, as desired.

$\textrm{(c)}\Rightarrow \textrm{(d)}$.
%Denote by
%$$
%P [\nu] (z) := \int_{\partial \mathbb D}
%P_{z}(\zeta) d\nu(\zeta) .
%$$
Following the original ideas of Beurling, and
similarly to the proof of Lemma 1 in \cite{HaLy}, for a given set of
indices $\mathcal Q\subset\mathcal I$ satisfying \eqref{caixes}
and for $(n,k)\in\mathcal Q$,
let $\lambda_{n,k} \in \Lambda\cap Q_{n,k}$, and define the (possibly divergent) series
$$
H_\nu (z) := \sum_{(n,k) \in \mathcal Q} 2^{-n}
P [\nu ] (\lambda_{n,k})
P_z(\lambda_{n,k}^* ),
$$
where $\lambda_{n,k}^*=\lambda_{n,k}/|\lambda_{n,k}|$.

The terms of this series are positive harmonic functions so by Harnack's theorem $H_\nu$ is either identically $+\infty$ or it defines a positive harmonic function.

Suppose that $\mathcal Q$ is a set for which
\eqref{HLint} fails. Then the series defining
$H_\nu(0)$ converges, and therefore $H_\nu$
is a positive harmonic function. Notice also that
$\lim_{|z|\to 1}(1-|z|)H_\nu(z)=0$,
since each term of the sum has this property,
and we can apply dominated convergence.

By retaining only the $(n,k)$ term of the sum, we see that
$$
H_\nu (\lambda_{n,k})
\succeq 2^{-n}  P [\nu] (\lambda_{n,k}) \frac1{1-|\lambda_{n,k}|}
  \succeq  P [\nu] (\lambda_{n,k}).
$$

Thus, using Harnack's inequality and  choosing an appropriate constant $C_0>0$,
we obtain a function
$$
h_\nu := P [\nu ] - C_0 H_\nu
$$
which is non-positive on $\bigcup_{(n,k)\in\mathcal Q}
Q_{n,k}$, and tends to infinity as $z$ tends radially
to the boundary a.e. with respect to $\nu$.

Now suppose that (d) doesn't hold.
This means that we are given a Blaschke
distribution $\{N_{n,k}\}_{(n,k)\in\mathcal I}$ and a singular measure
$\nu$ such that
for any $\mathcal Q$ such that \eqref{invul} holds, then \eqref{HLint} fails.

\textit{Claim.}
There exist a constant $\gamma >0$ and a
subset $\Lambda_0 \subset \Lambda$ such that
\begin{align*}
&\textrm{(i)} \qquad &&\sum_{(n,k) : (\Lambda \setminus \Lambda_0) \cap
Q_{n,k} \neq \emptyset} 2^{-n} P[\nu] (\lambda_{n,k}) 
\simeq \sum_{(n,k) : (\Lambda \setminus \Lambda_0) \cap
Q_{n,k} \neq \emptyset} 2^{-n} \int_{\T}  P_{c_{n,k}}(\zeta) d\nu(\zeta)
< \infty ,\\
&\textrm{(ii)} \qquad && \gamma P [\nu ] (\lambda_{n,k})
\le w_{n,k}(\Lambda, N_{n,k}) \ 
\textrm{for any $(n,k)$ with $\Lambda_0 \cap Q_{n,k} \neq \emptyset$}.
\end{align*}

To see this, define
$$
\mathcal L_j := \{ (n,k) \in \mathcal I : \Lambda
\cap Q_{n,k} \neq \emptyset \mbox{ and }
w_{n,k}(\Lambda , N_{n,k})  \le
2^{-j}  P [\nu ] (\lambda_{n,k}) \}.
$$
If there exists some $j_0$ such that
$$
\sum_{(n,k) \in \mathcal L_{j_0}} 2^{-n} 
P [\nu ] (\lambda_{n,k}) < \infty,
$$
then define $\Lambda \setminus \Lambda_0 := \Lambda
\cap \bigcup_{(n,k) \in \mathcal L_{j_0}}
Q_{n,k}$, and we have the result with $\gamma =
2^{-j_0} $.

Otherwise, set $j_1=1$ and define recursively
  $j_{m+1}>j_m$ and subsets $\mathcal A_m \subset \mathcal L_{j_{m}} $
such that
\begin{itemize}
\item
$\mathcal A_{m+1} \cap \mathcal A_l =\emptyset $, $1\le l \le m$,
and
\item
$1 \le \sum_{(n,k) \in \mathcal A_m} 2^{-n}
 P [\nu] (\lambda_{n,k}) \le M$.
\end{itemize}
This is possible because the terms to be summed
belong to a divergent series, and are all bounded
by a constant. Now, taking $\mathcal Q:= \cup
\mathcal A_m$, we have
$$
\sum_{(n,k) \in \mathcal Q} 2^{-n} w_{n,k}(\Lambda , N_{n,k})
\le \sum_m 2^{-j_m} \bigl(\sum_{(n,k) \in \mathcal A_m} 2^{-n}
P [\nu ] (\lambda_{n,k})\bigr)\leq
M \sum_m 2^{-j_m} < \infty,
$$
  while $ \sum_{(n,k) \in \mathcal Q} 2^{-n}
P [\nu] (\lambda_{n,k}) = \infty$,
which contradicts the hypothesis. The claim is
proved.

We now proceed to prove that the set $\Lambda$ is
not of determination for the Nevalinna class.  Let $h_\nu$ be
the function constructed above, using the set
$\Lambda \setminus \Lambda_0$ as the set which
doesn't satisfy \eqref{HLint}. Let $B$ be a Blaschke product with $N_{n,k}$ zeros $b_1, \dots, b_{N_{n,k}}$
located in $Q_{n,k}^\delta$ chosen as the
solution to the extremal problem in the
definition of $w_{n,k}$:
$$
\inf_{\lambda \in \Lambda \cap \overline Q_{n,k}} \sum_{j=1}^{N_{n,k}}
\log\frac1{\rho(\lambda, b_j)}
=
\sup_{a_1, \dots, a_{N_{n,k}} \in Q_{n,k}^\delta}
\inf_{\lambda \in \Lambda \cap \overline Q_{n,k}} \sum_{j=1}^{N_{n,k}}
\log\frac1{\rho(\lambda, a_j)}
=w_{n,k}(\Lambda , N_{n,k}) .
$$
Choose an integer $m$ such that $m\gamma \ge
1$ and pick a function $f \in \NN$ with
$$
\log |f| = m \log |B| + h_\nu .
$$
By construction
$h_\nu \le 0$ on $\Lambda\setminus\Lambda_0$
and $h_\nu \le \mathcal P[\nu]$ on $\Lambda$. Also, (ii) implies that
$\log|f|\leq 0$ on $\Lambda_0$, so altogether $\log
|f| \le 0$ on the whole of $\Lambda$. On the other hand the
fact that $\limsup_{|z|\to 1} (1-|z|) [\log |B(z)|+H_\nu(z)]= 0$
shows that $f$
cannot be bounded on the disk.

$\textrm{(d)}\Rightarrow \textrm{(e)}$. 
Condition (e) is the special case
of (d) where $\nu$ is a point mass.
\end{proof}

\section{Regular sampling sequences} \label{S3}

In this section we give precise conditions for three types of regular sets
to be sampling for $\NN$.

\subsection{Fine nets}\label{finenets}
Let $g:(0,1]\longrightarrow(0,1]$ be a non-decreasing continuous function with $g(0)=0$
A sequence $\Lambda$ is called a $g$\emph{-net} if
and only if 
\begin{itemize}
\item[(i)] The disks
$D (\lambda,   g(1-|\lambda|))$, $\lambda\in\Lambda$, are mutually disjoint,
\item[(ii)] There exists $C>0$ such that 
$\bigcup_{\lambda\in\Lambda} D (\lambda, C  g(1-|\lambda|)) = \D$.
\end{itemize}

We characterize sampling $g$-nets in terms of the growth of
$g$ (Theorem~\ref{thmxarxes} below).The
growth condition is equivalent to a condition in terms of approach regions.
Let 
\[
\mathcal F=\bigl\{
\psi:[0,1)\longrightarrow \R_+  \ \textrm{non-decreasing, continuous,  
with}\ \psi(0)=0\ \textrm{and}\  \int_0 \psi(x)/x^2\; 
dx<\infty\bigr\}.
\]
Given $\zeta\in\T$ define the approach region
$\Gamma_\psi(\zeta)=\{z\in\D : \psi(|z-\zeta|)\leq 1-|z| \}$.

\begin{theorem}\label{thmxarxes}
Let $\Lambda$ be a $g$-net.
The following properties are equivalent:
\begin{itemize}
\item[(a)] $\Lambda$ is a sampling sequence for $\NN$
\item[(b)] $\sum_{\lambda\in\Lambda\cap\Gamma_\psi(\zeta)} 1-|\lambda|=\infty$ 
for all $\zeta\in\T$ and all $\psi\in\mathcal F$.
\item[(c)] $\displaystyle \int_0 \frac{dt}{t^{1/2} g(t)} = \infty$.
\end{itemize}
\end{theorem}

\begin{remark}\label{Mn} (i) The conditions above can be reformulated in
terms of the number  $M_{n,k}$ of points of $\Lambda$ in a dyadic square
$Q_{n,k}$.   In this case $M_{n,k}$  is essentially independent of $k$, in
the sense that there exist $M_n$ and a constant $C>0$ such that $C^{-1}
M_n\leq M_{n,k}\leq C M_n$ for all  $0\leq k<2^n$.
Then the conditions in the theorem above are equivalent to $\sum_n (M_n
2^{-n})^{1/2}=\infty$ (see Lemma~\ref{Mnk}). 

(ii) As we will see soon, condition (b)  is always necessary. 
However, the results of Section \ref{uniformly} show that it is not always
sufficient.
\end{remark}

We will begin by proving that condition (b) is necessary in general. We first
give a reformulation of Example 6.4 in \cite{HaLy}, where we
look at when the outside of an approach region
$\Gamma_\psi(\zeta)$  is, or is not, Hayman-Lyons. All computations should be
done in the disc, depending only on what happens in a neighbourhood of a point
$\zeta$ on the boundary.  Passing to  the upper half-plane $U$ with the standard
conformal mapping, we may perform the corresponding computations in a disc of
fixed radius centered on any point of the real axis.

\begin{lemma}\label{regiotgt}
Let  $\psi:(0,\infty)\longrightarrow (0,1]$ be a non-decreasing 
continuous function with $\psi(0)=0$. 
%and  $\sup\limits_{x\in (0,\infty)} \eta (2x)/ \eta(x)<\infty$.
The set $\mathcal D_\psi := \{ x+iy \in \C : 0<y< \psi(|x|)
\mbox{ or } y \ge 1\}$ is a Hayman-Lyons set if and only if
$$
\int_0 \frac{\psi (x)}{x^2} \, dx = \infty.
$$
Furthermore, if $\mathcal D_\psi$ is not a Hayman-Lyons set, there
exists a harmonic function $h\in\Har_{\pm}(U)$, non-positive on  $\mathcal
D_\psi$, and with $\liminf\limits_{y\to 0} y h(iy)>0$.
\end{lemma}

\begin{proof}
For any point $\zeta \in \mathbb R$ except the origin (but including 
the point at
infinity), $\mathcal D_\psi$ contains a half-disc centered at $\zeta$, so
that the integral $I(\mathcal D_\psi^\delta,\zeta, \delta)\geq 
I(\mathcal D_\psi,\zeta, \delta) = \infty$. There
remains the case $\zeta=0$. 

A direct computation shows that
$$
I(\mathcal
D_\psi,0, 0) =\int_{\mathcal D_\psi}\frac{dx\, dy}{x^2+y^2}\simeq
\int_0^1  \int_0^{\psi(x)}\frac1{x^2+y^2} =
\int_0^1 \frac{\arctan(\psi(x)/x)}x\; dx\ ,
$$
which is finite if and only if $\int_0  \psi (x)/x^2 \, dx<\infty$.
Hence $  I(\mathcal D_\psi,0, 0) =\infty$ when
$\int_0\psi(x)/x^2\; dx=\infty$.

This same estimate shows that in order to prove that $ I(\mathcal D_\psi^\delta,0, 0) <\infty$ when
$\int_0\psi(x)/x^2\; dx<\infty$ it is enough to see 
that for any $\delta>0$, there exists a non-decreasing  function
$\psi_\delta\geq 0$ with $\psi_\delta(0)=0$ and such that
\begin{itemize}
\item[(i)] $\De^\delta\subset\mathcal D_{\psi_\delta}$
\item[(ii)] $\int_0 \psi_\delta(x)/x^2\; dx\simeq \int_0 \psi(x)/x^2\; dx$.
\end{itemize}

In the construction of $\psi_\delta$ only the behavior near zero is relevant, hence we restrict
our attention to $x\in[0,1/2]$. Recall that 
the pseudohyperbolic distance between two points $z,\zeta\in U$ 
is given by $\rho(z,\zeta)=|z-\zeta|/|z-\bar\zeta|$.

Let $\eta>0$ (to be chosen later) and consider the function
\[ 
\psi_\eta(x)=\sum_{n=0}^\infty \frac{1+\eta}{1-\eta}
\psi(2^{-n})\chi_{[2^{-(n+2)},2^{-(n+1)})}(x)\ .
\]
This corresponds to ``raising" the value of $\psi$ at  $x=2^{-n}$ by $\eta$ 
(in the pseudohyperbolic metric) and assign it to the whole interval 
$[2^{-(n+2)},2^{-(n+1)})$ (see Figure 2).

It is clear that $\psi_\eta$ satisfies (ii) for any $\eta>0$:
\[
\int_0 \frac{\psi_\eta(x)}{x^2} dx= \frac{1+\eta}{1-\eta}
\sum_{n=0}^\infty \psi(2^{-n})\int_{2^{-(n+2)}}^{2^{-(n+1)}}\frac{dx}{x^2}
\simeq \sum_{n=0}^\infty \frac{\psi(2^{-n})}{2^{-n}}\simeq \int_0 \frac{\psi(x)}{x^2} dx\ .
\]

Let us see that the pseudohyperbolic distance between the graph of 
$\psi_\eta$ and the graph of $\psi$ is bigger than $\delta$ if $\eta$ is big enough, and therefore (i) holds as well. In the vertical direction it is clear that we only need to take $\eta\geq\delta$, by construction of $\psi_\eta$. For the horizontal direction we have, for 
any $x \in [2^{-(n+1)}, 2^{-n})$, 
\begin{align*}
\rho\bigl((2^{-(n+2)},\psi_\eta(2^{-(n+2)}- )), (x,\psi(x))
\bigr)
&=\frac {\left|x-2^{-(n+2)} + i (\psi(x)-\frac{1+\eta}{1-\eta}\psi(2^{-(n+1)}))\right|}{\left|x-2^{-(n+2)}+ i(\psi(x)+\frac{1+\eta}{1-\eta}\psi(2^{-(n+1)}))\right|}\\
&\ge
\frac{1}{\left|1 + i\ \frac{1+\eta}{1-\eta} \frac{\psi(x)+ \psi(2^{-(n+1)})}{x-2^{-(n+2)}}\right|}\ 
\end{align*}
This is clearly bounded below, since the integrability condition on $\psi$ gives
in particular 
\[
\lim_{n\to\infty}\frac{\psi(x)+ \psi(2^{-(n+1)})}{x-2^{-(n+2)}} \le
\lim_{n\to\infty}\frac{2 \psi (2^{-n})}{2^{-(n+1)}-2^{-(n+2)}}=
\lim_{n\to\infty}\frac{8 \psi(2^{-n})}{2^{-n}}=0\ .
\]

\vglue 0.5 truecm
\begin{center}
\psfrag{tag0}{\scriptsize$0$}
\psfrag{tag1}{\tiny$2^{-n}$}
\psfrag{tag2}{\tiny$2^{-(n+1)}$}
\psfrag{tag3}{\tiny$2^{-(n+2)}$}
\psfrag{tag4}{\tiny$2^{-(n+3)}$}
\psfrag{tag5}{\scriptsize{$\psi(|x|)$}}
\psfrag{tag6}{\scriptsize{$\mathcal D_\psi$}}
\psfrag{tag7}{\scriptsize{$\eta$}}
\psfrag{tag8}{\scriptsize{$\psi_\eta(|x|)$}}
\psfrag{tag9}{\scriptsize{$\mathbb R$}}
%\psfrag{...}{$\cdots$}
\includegraphics[totalheight=5.5cm]{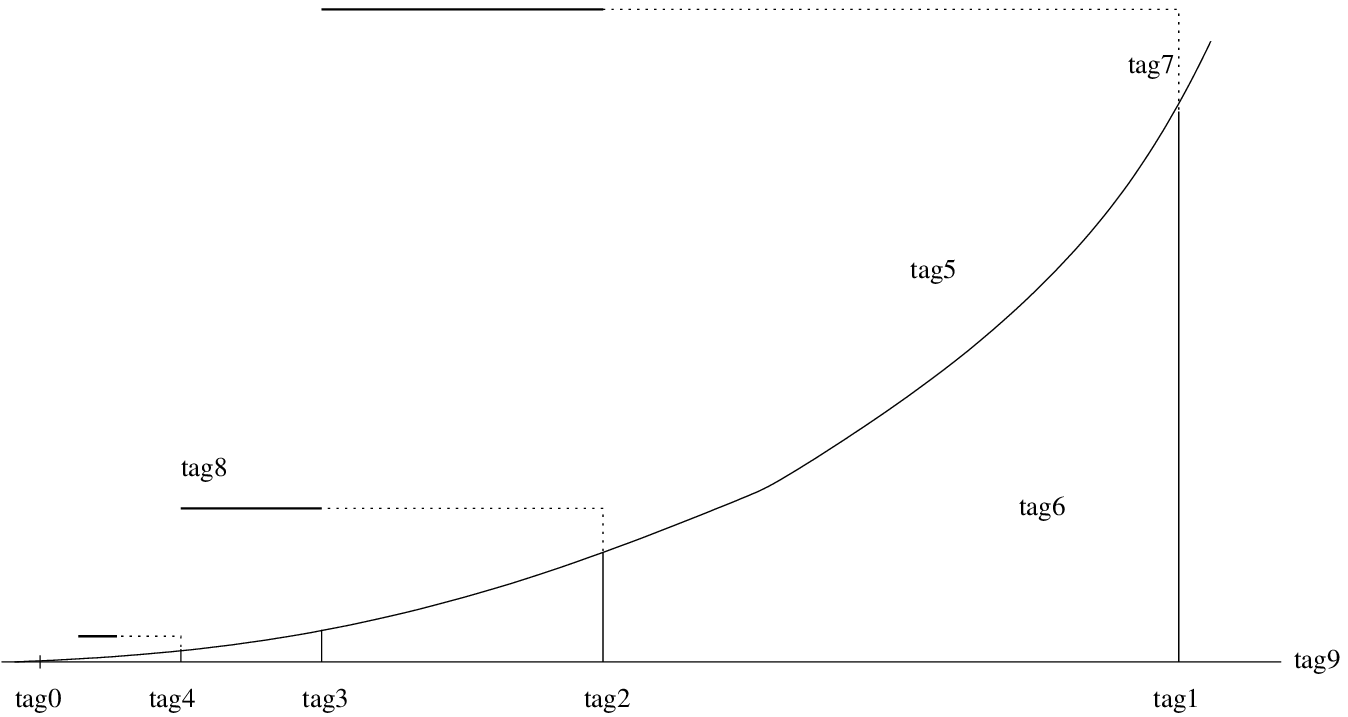}

{\footnotesize{Figure 2.}}
%{\footnotesize{Figure 1. Upper-half plane representation of $\Lambda$ and the horocycle $P_z(1)=c$.}}
\end{center}
\vglue 0.5 truecm

In the case where the integral is convergent, denoting temporarily
$P_{x+iy}(t) := \frac1\pi \frac{y}{(x-t)^2+y^2}$, let
$$
h (x+iy) := P_{x+iy}(0) - C h_\psi (x+iy),
$$
where $C>0$ and $h_\psi$ is the Poisson integral of the integrable function $\psi(t)/t^2$
restricted to the interval $[-1,1]$. 

It is well known that the
growth of the Poisson integral of an integrable function is such that
$\lim_{y\to 0} y h_\psi(x+iy)=0$ (see for instance
\cite[Corollary 2]{SS}), thus it will be enough to prove that for $C>0$ sufficiently big
$h\leq 0$ on $\mathcal D_\psi$.

This will be done as soon as we see that $h\leq 0$ for
points $z=x_0+iy_0\in\partial \mathcal D_\psi$. Since $\lim\limits_{x\to 0}\psi(x)/x=0$, we have then
\[
P_z(0)=\frac{\psi(x_0)}{x_0^2+\psi^2(x_0)}\simeq\frac{\psi(x_0)}{x_0^2}\ .
\]
On the other hand
\begin{align*}
h_\psi(x_0+iy_0)&\succeq \int_{|x-x_0|<\psi(x_0)}\frac 1{\psi(x_0)}\frac{\psi(x)}{x^2}\, dx
\succeq \frac 1{\psi(x_0)} \int_{x_0}^{x_0+\psi(x_0)} \frac{\psi(x)}{x^2}\, dx\\
&\geq \bigl(\frac 1{x_0}-\frac 1{x_0+\psi(x_0)}\bigr)=\frac{\psi(x_0)}{x_0^2}\frac 1{1+\psi(x_0)/x_0}\simeq
\frac{\psi(x_0)}{x_0^2}\ ,
\end{align*}
hence with $C$ big enough we get the desired estimate.

\end{proof}

\begin{proof}[Proof of Theorem~\ref{thmxarxes}]
(a)$\Rightarrow$(b). 
Assume that there exist a function $\psi\in \mathcal F$ and $\zeta\in\T$ such that
$\sum_{\lambda\in\Lambda\cap\Gamma_\psi(\zeta)} 1-|\lambda|<\infty$.
Similarly to Example 1, consider the function $f\in\NN$ such that
\[
\log |f|=\log|B_Z|+h_\psi\ ,
\]
where $B_Z$ is the Blaschke product associated to $Z=\Lambda\cap
\Gamma_\psi(\zeta)$ and $h_\psi$ is the harmonic function obtained by
transporting to the disk the function given by the previous
lemma. 

It is clear then that $f$ is bounded on $\Lambda$. But by the properties
of $h_\psi$ and by Lemma~\ref{blaschke} we see that $f$ cannot be bounded in the whole disk.
Thus $\Lambda$ is not a determination set for $\NN$.

(b)$\Rightarrow$(c). It will be enough to prove the following lemma.

\begin{lemma}\label{Mnk}
Assume that $\Lambda$ satisfies condition (b) in 
Theorem~\ref{thmxarxes}. Let $M_{n,k}=\#\Lambda\cap Q_{n,k}$
and assume that 
there exist $C>0$ and positive integers $M_n$ such that
$C^{-1} M_n\leq M_{n,k}\leq C M_{n}$. 
Then $\sum_n(M_n 2^{-n})^{1/2}=\infty$.
\end{lemma}

A volume estimate shows that
$M_n\simeq (g(2^{-n}))^{-2}$, and therefore
\[
\int_0 \frac{dt}{t^{1/2} g(t)}\simeq \sum_n 2^{-n/2}\frac 1{g(2^{-n})}\simeq
\sum_n(M_n 2^{-n})^{1/2}=\infty\ .
\]

\begin{proof}
We want to prove that if $\sum_n(M_n 2^{-n})^{1/2}
<\infty$ there exists $\psi\in\mathcal F$ such that
$\sum_{\lambda\in\lambda\cap \Gamma_\psi(\zeta)} 1-|\lambda|<\infty$, thus
contradicting (b).

It will be expedient to write
the computation in the upper half-plane $U$ after a conformal mapping. 
We call the resulting sets $\Lambda$ and $\Gamma_\psi(\zeta)$ again. 
Thus, we need to prove that there exists $\psi\in\mathcal F$
\begin{equation}\label{suma}
\sum_{\lambda\in\Lambda\cap\Gamma_\psi(\zeta)} 1-|\lambda|\simeq
\sum_{n\ge 0} 2^{-n} \#
\left(
\Lambda \cap \Gamma_\psi(\zeta) \cap \{2^{-n-1} < y \le 2^{-n} \}
\right)<\infty\ .
\end{equation}

Since $\psi$ is an increasing function,  the set
$\Gamma_\psi(\zeta) \cap \{2^{-n-1} < y \le 2^{-n} \}$ is contained in the rectangle
$\{ |x| \le \psi^{-1}(2^{-n}), 2^{-n-1} < y \le 2^{-n} \}$.
Therefore, splitting the sum for the different $Q_{n,k}$ and using that 
$M_n\simeq g^{-2}(2^{-n})$
we have
\begin{align*}
\sum_{\lambda\in\Lambda\cap\Gamma_\eta(\zeta)} 1-|\lambda|&\preceq \sum_{n\geq 0}2^{-n}
\sum_{k: k 2^{-n}\leq\psi^{-1}(2^{-n})} M_{n,k}
%\preceq \sum_{n\geq 0}2^{-n}
%M_n\frac{\Psi^{-1}(2^{-n})}{2^{-n}}\\
\preceq \sum_{n\geq 0}  M_n \psi^{-1}(2^{-n})
\simeq \sum_{n\geq 0}  \frac{\psi^{-1}(2^{-n})}{g^2(2^{-n})}
\ .
\end{align*}

It will be enough to see that $\psi$
defined by $\psi^{-1}(t)=\sqrt t g(t)$ is in $\mathcal F$, since then
\[
\sum_{\lambda\in\Gamma_\eta(\zeta)} 1-|\lambda|\preceq 
\sum_{n\geq 0} \frac{2^{-n/2}}{g(2^{-n})}\simeq \int_0 \frac{dt}{\sqrt t g(t)}<\infty\ .
\]
By definition $\psi$ is non-decreasing, continuous and $\psi(0)=0$. Also
\begin{align*}
\int_0^1 \frac{dt}{\sqrt t g(t)}&=\int_0^1 \frac {dt}{\psi^{-1}(t)}
=\int_1^\infty \left|\{ t : 1/\psi^{-1}(t)\geq \alpha\}\right|\, d\alpha 
=\int_0^1 \left|\{ t : \psi^{-1}(t)\leq s\}\right|\, \frac{ds}{s^2}\\
&=
\int_0^1 \left|\{ t : t\leq \psi(s)\}\right|\, \frac{ds}{s^2}
= \int_0^1 \frac{\psi(s)}{s^2}\; ds<\infty\ ,
\end{align*}
as desired.
\end{proof}

(c)$\Rightarrow$(a).  Given a Blaschke sequence $Z$ we want to choose a
family of indices $\mathcal Q\subset\mathcal I$
%$\{Q_{n,k}\}_{(n,k)\in\mathcal Q}$ 
satisfying the conditions of Theorem
\ref{suficient}(a). First we need a control of the vulnerability on the squares where 
$N_{n,k}=\#(Z\cap\tqnk)$ is small.

\begin{lemma}\label{west}
If $\Lambda$ is a $g$-net,  there exist $\eps >0 $, $n_0 \in \N$ and  $C > 0 $
such that whenever $n \ge n_0$ and
$N \le \eps \# (\Lambda \cap Q_{n,k})$, then
$w_{n,k}(\Lambda , N ) \le C N$.
\end{lemma}

\begin{proof}
Denote $Q=Q_{n,k}$, $\tilde Q=\tilde Q_{n,k}$ and
$M=\# (\Lambda \cap Q)$. Then
$M\simeq g(2^{-n})^{-2}$ and 
there exist constants $c < C$ such that the disks
$D(\la,c g(2^{-n}))$, $\lambda\in Q$, are mutually disjoint and the disks
$D(\la,C g(2^{-n}))$ cover the whole set $Q$.
Let $Z=\{a_1,\dots,a_N\}\subset \overline{Q}^\delta$ and consider
\[
\Lambda_Q^\prime:= \bigl\{ \la \in \Lambda \cap Q : \rho (\lambda, a) 
\ge C g(2^{-n}),
\mbox{ for all } a \in Z \cap \tilde Q \bigr\}.
\]
For $\eps$ small enough, and $n$ (and therefore $M$) large enough,
$\# \Lambda_Q^\prime \ge M/2$.

Choose $r\in (0,1)$  large enough so that $\overline{Q} \subset D(a, r)$ 
for any $a \in \tilde Q$.
For $\lambda'\in \Lambda_Q^\prime$ and $z \in D(\la', c
g(2^{-n}))$ we have $\rho(z, a_j) \preceq \rho(\lambda', a_j)$ for all $a_j \in \tilde Q
$, and therefore 
\begin{equation}\label{average}
\log \frac1{\rho(\la', a_j)} \preceq \frac 1{m(D(\la', c g(2^{-n}))}
\int_{D(\la', c g(2^{-n}))}\log \frac1{\rho(z, a_j)}\; dm(z).
\end{equation}
Then, for any $a_1, \dots, a_N \in \tilde Q $,
\begin{align*}
\sum_{\la' \in \Lambda_Q^\prime} \sum_{j=1}^N \log \frac1{\rho(\la', a_j)}
=
\sum_{j=1}^N \sum_{\la' \in \Lambda_Q^\prime} \log \frac1{\rho(\la', a_j)}
\preceq
\sum_{j=1}^N  \frac1{( g(2^{-n}) 2^{-n})^2} \int_{D(a_j, r)}
\log \frac1{\rho(z, a_j)} dm(z),
\end{align*}

Applying an automorphism of the disk sending $a_j$ to
the origin, we see that there is a radius $r'\in (r,1)$ with
$$
\frac1{(2^{-n})^2}
\int_{D(a_j, r)}
\log \frac1{\rho(z, a_j)} dm(z)
\preceq
\int_{D(0, r')}
\log \frac1{|z|} dm(z) \le C.
$$
Finally,
$$
\sum_{\la' \in \Lambda_Q^\prime} \left( \sum_{j=1}^N \log 
\frac1{\rho(\la', a_j)}  \right)
\preceq
\frac 1{g^2(2^{-n})} N \simeq M N,
$$
and since $\# \Lambda_Q^\prime \ge M/2$, the average value of the summands in the
first sum is bounded by a constant multiple of $N$. 
\end{proof}

For $\eps\in(0,1)$ small enough, define
\begin{equation}\label{Q}
\mathcal Q=\{(n,k) : N_{n,k}\leq \eps \# (\Lambda \cap Q_{n,k})\}\ .
\end{equation}
Lemma \ref{west} together with the Blaschke condition give \eqref{invul}
in Theorem~\ref{suficient}.

In order to prove \eqref{HLzeta} in Theorem~\ref{suficient},
for each $(n,k)$ pick $c_{n,k} \in Q_{n,k}$, for instance its center.
By rotation invariance of the properties of being a $g$-net or a Blaschke
sequence, it will be enough to see that,
\begin{equation} \label{HLRbe}
\sum_{(n,k) \in \mathcal Q} \left(
\frac{1-|c_{n,k}|^2}{|1-c_{n,k}|} \right)^2 = \infty.
\end{equation}
Let
\begin{equation}\label{Ln}
L_n := \# \{ k : (n,k) \notin \mathcal Q \},
\end{equation}
Observe that $\frac{1-|c|^2}{|1-c|}$ is  bigger when the argument of $c$
in $(-\pi,\pi]$
is closer to $0$.
So we must have, for any fixed $n$,
\begin{equation}\label{sumLn}
\sum_{k : (n,k) \in \mathcal Q} \left(
\frac{1-|c_{n,k}|^2}{|1-c_{n,k}|} \right)^2
\succeq
\sum_{k\ge L_n}
\left(
\frac{2^{-n}}{k2^{-n}} \right)^2
\simeq \frac1{L_n}.
\end{equation}
On the other hand
\[
\sum_{k : (n,k)\notin \mathcal Q} N_{n,k}\geq 
\eps \sum_{k : (n,k)\notin \mathcal Q} \#(\Lambda\cap Q_{n,k})\geq
\eps C \frac{L_n }{(g(2^{-n}))^2},
\]
and the Blaschke condition implies
$$
\sum_n \frac{ 2^{-n} L_n}{g(2^{-n})^2} < \infty.
$$

If \eqref{HLRbe} doesn't hold we have
$\sum_n 1/L_n<\infty$, and
$$
\int_0\frac{dt}{t^{1/2} g(t)}\simeq
\sum_{n\geq 0}\frac{ 2^{-n/2} }{g(2^{-n})} =
\sum_{n\geq 0}
\left(
\frac1{L_n}
\frac{ 2^{-n} L_n}{g(2^{-n})^2}
\right)^{1/2}\leq \bigl(\sum_{n\geq 0}\frac1{L_n}\bigr)^{1/2}
\bigl(\sum_{n\geq 0} \frac{ 2^{-n} L_n}{g(2^{-n})^2}\bigr)^{1/2}
< \infty, 
$$
which contradicts the hypothesis.
\end{proof}

\subsection{Discretized rings}
\label{d-rings}
Let $r_n\in (0,1)$ be an increasing sequence of radii with $\lim_n r_n=1$ and
$0<\inf_n \frac{1-r_{n+1}}{1-r_n}<\sup_n \frac{1-r_{n+1}}{1-r_n}<1$. Let $\epsilon_n$ be a decreasing sequence
of hyperbolic distances such that $\lim_n \epsilon_n=0$. The 
\emph{discretized rings associated
to $\{r_n\}_n$ and $\{\epsilon_n\}_n$} is the sequence $\Lambda=\{\lambda_{n,j}\}_{n,j}$,
where
\[
\lambda_{n,j}=r_n \exp\left( j \frac{2\pi i}{(1-r_n)\epsilon_n} \right)\qquad n\in\N,
\quad 0\leq j<\left[\frac{1}{(1-r_n)\epsilon_n}\right]\ .
\]

\begin{theorem}\label{thmlines}
Let $\Lambda$ be discretized rings. The following properties are equivalent:
\begin{itemize}
\item[(a)] $\Lambda$ is a sampling sequence for $\NN$
\item[(b)] $\sum_{\lambda\in\Lambda\cap\Gamma_\eta(\zeta)} 1-|\lambda|=\infty$ 
for all $\zeta\in\T$ and all $\psi\in\mathcal F$.
\item[(c)] $\displaystyle \sum_{n=0}^\infty\left(\frac{1-r_n}{\epsilon_n}\right)^{1/2}=\infty$.
\end{itemize}
\end{theorem}

The proof follows the same scheme as the proof of Theorem~\ref{thmxarxes}.

\begin{proof} 
(a)$\Rightarrow$(b). As in Theorem~\ref{thmxarxes}.

(b)$\Rightarrow$(c). Assume that $\sum_n\left(\frac{1-r_n}{\epsilon_n}\right)^{1/2}<\infty$.
Consider the  sequence $\eta_n=[(1-r_n)\epsilon_n]^{1/2}$ decreasing to zero and the function $\psi$
which in each interval $[\eta_{n+1},\eta_n)$ is defined as the segment joining
the points $(\eta_{n+1},1-r_{n+1})$ and $(\eta_{n},1-r_{n})$; i.e.
\[
\psi(x)=\sum_{n\geq 0} \bigl[(1-r_{n})+\frac{r_{n+1}-r_{n}}{\eta_{n}-\eta_{n+1}}(x-\eta_n)\bigr]\chi_{[\eta_{n+1},\eta_n)}(x)\ .
\]
It is clear that $\psi$ is continuous, non-decreasing and with $\psi(0)=0$. Also,
by the hypotheses on $\{r_n\}_n$
\begin{align*}
\int_0\frac{\psi(x)}{x^2}dx &= \sum_{n\geq 0}\int_{\eta_{n+1}}^{\eta_n} \frac{\psi(x)}{x^2} dx\leq
\sum_{n\geq 0} (1-r_n) \int_{\eta_{n+1}}^{\eta_n} \frac{dx}{x^2}
\leq \sum_{n\geq 0}(1-r_{n})\left(\frac 1{\eta_{n+1}}-\frac
1{\eta_n}\right)\\
&\leq \sum_{n\geq
0}\frac{1-r_{n+1}}{\eta_{n+1}}=\sum_n\left(\frac{1-r_n}{\epsilon_n}\right)^{1/2}<\infty ,
\end{align*}
hence $\psi\in\mathcal F$.

On the other hand, a length estimate shows that the number of $j$ such that
$\lambda_{n,j}\in\Gamma_\psi(\zeta)$ is approximately
$\frac{\psi^{-1}(1-r_n)}{(1-r_n)\epsilon_n}$. Since $\psi^{-1}(1-r_n)=\eta_n$ we have
\[
\sum_{\lambda\in\Lambda\cap\Gamma_\psi(\zeta)} 1-|\lambda|\simeq \sum_{n\geq 0}(1-r_{n}) 
\frac{\psi^{-1}(1-r_n)}{(1-r_n)\epsilon_n}=\sum_{n\geq 0}\frac{\eta_n}{\epsilon_n}=
\sum_{n\geq 0} \left(\frac{1-r_n}{\epsilon_n}\right)^{1/2}<\infty\ ,
\]
which contradicts (b).

(c)$\Rightarrow$(a). Given a Blaschke sequence $Z$, we want to choose a family of indices
$\mathcal Q$ satisfying the conditions of Theorem~\ref{suficient}(a). 

We begin by showing that Lemma~\ref{west} still holds for discretized rings.

\begin{proof}[Proof of Lemma~\ref{west} for discretized rings]
Let $Q=Q_{n,k}$, $\tilde Q=\tilde Q_{n,k}$ and $M=\# (\Lambda \cap Q)$. 

Since $0<\inf_n\frac{1-r_{n+1}}{1-r_n}<\sup_n\frac{1-r_{n+1}}{1-r_n}<1$, there is at most a finite number of $r_m$ with
$2^{-(n+1)}<1-r_m\leq 2^{-n}$. There is no loss of generality in assuming that there is just
one such $r_m$, and therefore $M\simeq 1/\epsilon_m$. In particular, 
there exist  $c < C$ such that  $D(\la,c \epsilon_m )$, $\lambda\in Q$, are
mutually disjoint and  $D(\la,C \epsilon_m)$ cover the segment 
$ \{z : |z|=r_m\}\cap Q$. 

Given $Z=\{a_1,\dots,a_N\}\subset \overline{Q}^\delta$ consider now
\[
\Lambda_Q^\prime:= \bigl\{ \la \in \Lambda \cap Q : \rho (\lambda, a) 
\ge C \epsilon_m,
\mbox{ for all } a \in Z \cap \tilde Q \bigr\}.
\]
Again, for suitable $C$, $m$ and $r\in (0,1)$, we have $\# \Lambda_Q^\prime \ge M/2$
and $\overline{Q} \subset D(a, r)$  for any $a \in \tilde Q$.
We proceed as before, but replacing the area averages in \eqref{average} by the
line averages
\begin{align*}
\log \frac1{\rho(\la', a_j)} &\preceq\frac 1{|J(\lambda',c \epsilon_m)|}
\int_{J(\lambda',c \epsilon_m)}\log \frac1{\rho(z, a_j)}\; |dz|\\
&\simeq
\frac 1{(1-r_m)\epsilon_m} \int_{J(\lambda',c \epsilon_m)}\log \frac1{\rho(z, a_j)}\; |dz| ,
\end{align*}
where 
$J(\lambda',c \epsilon_m)=\{z : |z|=r_m\}\cap D(\lambda',c \epsilon_m)$. Then,
for any $a_1, \dots, a_N \in \tilde Q $,
\begin{align*}
\sum_{\la' \in \Lambda_Q^\prime} \sum_{j=1}^N \log \frac1{\rho(\la', a_j)}
&=
\sum_{j=1}^N \sum_{\la' \in \Lambda_Q^\prime} \log \frac 1{\rho(\la', a_j)}
%\sum_{j=1}^N  \frac1{ (1-r_m)\epsilon_m} \sum_{\la' \in \Lambda'} 
%\int_{J(\lambda',c \epsilon_m)}\log \frac1{\rho(z, a_j)}\; |dz|
\preceq \sum_{j=1}^N  \frac1{ (1-r_m)\epsilon_m} 
\int_{Q\cap\{|z|=r_m\}}
\log \frac1{\rho(z, a_j)} |dz|\ .
\end{align*}
Since the hyperbolic length of $Q\cap\{|z|=r_m\}$ is approximately $2^{-n}\simeq 1-r_m$, we
have, for some $r'<1$,
\[
\frac 1{1-r_m} \int_{Q\cap\{|z|=r_m\}}
\log \frac1{\rho(z, a_j)} |dz|\preceq
\int_{-r'}^{r'} \log\frac 1{|x|}\; dx\simeq 1\ .
\]
Hence
\[
\sum_{j=1}^N \sum_{\la' \in \Lambda_Q^\prime} \log \frac 1{\rho(\la', a_j)}
\preceq N/\epsilon_m\ ,
\]
and using that $\# \Lambda_Q^\prime\geq 1/(2\epsilon_m)$ we get the desired result.
\end{proof}

From here we proceed as in the proof of Theorem~\ref{thmxarxes}. 
As pointed out before, there is no loss of generality in assuming that there is just
one $r_m$ with $2^{n+1}<1-r_m\leq 2^{-n}$. To simplify the notation we re-index
$r_m$ and call it $r_n$.

Given a 
Blaschke sequence $Z$ and $\eps$ small enough, define $\mathcal Q$ as in \eqref{Q}.
The previous lemma ensures \eqref{invul} in Theorem~\ref{suficient}. In order to see
that \eqref{HLzeta} also holds it is enough to show \eqref{HLRbe}.
Let $L_n$ as in \eqref{Ln}. Since $\#(\Lambda\cap Q_{n,k})\simeq 1/\epsilon_n$ we have now
\[
\sum_{k : (n,k)\in\mathcal Q} N_{n,k}\geq\eps 
\sum_{k : (n,k)\in\mathcal Q} \#(\Lambda\cap Q_{n,k})\geq\eps CL_n/\epsilon_n,
\]
hence the Blaschke condition implies $\sum_n (1-r_n) L_n/\epsilon_{n}<\infty$\ .
If \eqref{HLRbe} does not hold we have $\sum_n 1/L_n<\infty$, and
\[
\sum_n \left(\frac{1-r_n}{\epsilon_n}\right)^{1/2}\leq \left(\sum_n\frac 1{L_n}\right)^{1/2}\left(\sum_n\frac{1-r_n}{\epsilon_n}L_n\right)^{1/2}<\infty\ ,
\]
which contradicts the hypothesis.
\end{proof}

\subsection{Uniformly dense disks}\label{uniformly}

In this section we consider a different kind of sampling sets. We begin 
with the sequences considered by Ortega-Cerd\`a and Seip in \cite{OrSe}.

\begin{definition}
A sequence $\Lambda\subset\D$ is \emph{uniformly dense} if
\begin{itemize}
\item[(i)] $\Lambda$ is separated, i.e. 
$\inf_{\lambda\neq \lambda'} \rho(\lambda,\lambda')>0$.
\item[(ii)] There exists $r<1$ such that $\D=\bigcup_{\lambda\in\Lambda} D(\lambda,r)$.
\end{itemize}

\end{definition}

Notice that, in the terminology of Section~\ref{finenets}, uniformly dense
sequences correspond to $1$-nets.

Let $\varphi$ be a non-decreasing continuous function, bounded by some
constant less than 1. Given $\Lambda$ uniformly dense define 
$r_\lambda=\varphi(1-|\lambda|)$, $D_\lambda^{\varphi}=D (\la , r_\lambda)$
and the unions of disks
$$
\Lambda(\varphi) := \bigcup_{\la\in\Lambda} D_{\la}^{\varphi} .
$$

\begin{theorem}\label{ud}
The set $\Lambda(\varphi)$ is sampling for $\NN$ if and only if
\begin{equation}\label{LyubSeip}
\int_0^1 \frac{dt}{t \log (1/\varphi(t))} =
\infty\ .
\end{equation}

\end{theorem}

In Section~\ref{S5} we will see that this condition actually characterizes
determination sets for the space of subharmonic functions in the disk having the 
characteristic growth of the Nevanlinna class.

\begin{remark}\label{sum}

Condition \eqref{LyubSeip} is equivalent to the fact that the harmonic
measure of $\partial\D$ in $\D \setminus \overline{\Lambda(\varphi)}$
is zero, see \cite[Theorem 1]{OrSe}. Notice also
that for any fixed $K>1$ condition \eqref{LyubSeip}
is equivalent to
\[
\sum_n \frac1{\log (1/\varphi(K^{-n}))}= \infty\ .
\]

\end{remark}

\begin{remark}\label{A-a}
The above family of examples allows
us to see that
there is no general relationship between $\aa$-sampling sets and Nevanlinna sampling sets.

A set $\Lambda\subset\D$ is \emph{sampling} for the space
\[
A^{-\alpha}=\{f\in\Hol (\D) : \|f\|_{\alpha}:=\sup_{z\in\D} (1-|z|)^\alpha|f(z)|<\infty\}\qquad \alpha>0,
\]
when there exists $C>0$ such that $\|f\|_{\alpha}\leq C 
\sup_{\lambda\in\Lambda} (1-|\lambda|)^\alpha|f(\lambda)|$
for all $f\in A^{-\alpha}$.

A well-known result of K. Seip \cite[Theorem 1.1]{Se} characterizes $\aa$-sampling sets as those $\Lambda$ for which there 
exists a separated subsequence $\Lambda'\subset\Lambda$ such that
\[
D_-(\Lambda'):=\liminf_{r\to 1^{-}}\inf_{z\in\D}
\frac{\sum\limits_{\lambda: 1/2<\rho(\lambda,z)<r}\log\frac 1{\rho(\lambda,z)}}{\log\frac 1{1-r}}>\alpha\ .
\]

Let $\Lambda_g$ be a fine net associated to a function $g$ with $\int_0\frac{dt}{t^{1/2} g(t)}<\infty$, for instance $g(t)=t^{1/4}$. According to Theorem~\ref{thmxarxes}, $\Lambda_g$ is not a Nevanlinna sampling set. On the other hand, for any given $\alpha>0$, we can extract a maximal separated sequence $\Lambda'$ with the separation small enough so that $D_-(\Lambda')>\alpha$, hence $\Lambda_g$ is $\aa$-sampling for all $\alpha>0$.

Also, given $\alpha>0$, consider a uniformly dense sequence $\Lambda$ with $D_-(\Lambda)<\alpha$
and take $\varphi$ satisfying $\lim_{t\to 0}\varphi(t)=0$ and \eqref{LyubSeip}. Then $\Lambda(\varphi)$ is Nevanlinna sampling but it is not $\aa$-sampling, since $D_-(\Lambda')<\alpha$ for any separated $\Lambda'\subset \Lambda(\varphi)$. 

Better yet, take a set $\Lambda$ as in Section \ref{d-rings}, sampling for the Nevanlinna class, with $\lim_{n\to \infty} \frac{1-r_n}{1-r_{n+1}} = \infty$. Then $D_-(\Lambda)=0$, so it cannot be $\aa$-sampling for any $\alpha>0$.

\end{remark}

\begin{proof}
Assume  that \eqref{LyubSeip} does not hold.
%$\int_0^1 \frac{dt}{t \log  (1/\varphi(t))}<\infty$.
We will exhibit a function $f\in \NN$ such that
$\log|f(z)|=  \log|B(z)|+\delta P_z(1)$  is bounded on $\Lambda(\varphi)$, for an
appropriate choice of the Blaschke product $B$ and the constant $\delta>0$. 
Since, according to Lemma~\ref{blaschke}, $f\notin \H^\infty$, this will contradict the
fact that $\Lambda(\varphi)$ is sampling.

Let $Z$ be the set of $\la\in\Lambda$ such that
$m_\la:=\bigl[
\frac{P_{\la}(1)}{\log  1/\varphi(1-|\la|)} \bigr]\geq 1$, 
where each point $\la$ is taken with multiplicity $m_{\la}$. 

In order to see that the Blaschke sum of $Z$ (with multiplicities)
is finite,  split it  into the different dyadic 
squares $Q_{n,k}$. Notice that for $\la\in Q_{n,k}$
\[
P_{\la}(1)=\frac{1-|\la|^2}{|1-\la|^2}\simeq
\frac{2^{-n}}{(2^{-n}+k2^{-n})^2}=\frac{2^n}{(1+k)^2}\ .
\]
Also, the uniform density condition implies
$\# \Lambda\cap Q_{n,k}\preceq 1$. Therefore
\begin{align*}
\sum_{a \in Z} (1-|a|) &\simeq
\sum_{n\ge 0} \sum_{k =0}^{2^n-1}\sum_{\la\in Z\cap Q_{n,k}}
m_\la (1-|\la|)
\simeq
\sum_{n\ge 0} 2^{-n}\sum_{k =0}^{2^n-1}\sum_{\la\in Z\cap Q_{n,k}}
\frac{ P_{\la}(1)}{\log (1/\varphi(1-|\la|))}\\
&\preceq
\sum_{n\ge 0} \frac{1}{\log (1/\varphi(2^{-n}))}
\bigl(\sum_{k =0}^{2^n-1}\frac1{1+k^2}\bigr)\simeq 
\int_0^1 \frac{dt}{t \log  (1/\varphi(t))}  <\infty.
\end{align*}
On the other hand, if $z\in D_{\la}^{\varphi}$ we have
\[
\log |B(z)|\leq \log[\rho(z,\la)]^{m_{\la}}\preceq
\frac{P_{\la}(1)}{\log  1/\varphi(1-|\la|)}\log \varphi(1-|\la|)=
-P_{\la}(1).
\]
Therefore $\log |f|$ is bounded on 
$\bigcup_{\la\in\Lambda} D_{\lambda}^{\varphi}$ if $\delta$ is chosen small enough.

Assume now that \eqref{LyubSeip} holds.
%$\int_0^1 \frac{dt}{t \log  (1/\varphi(t))}=\infty$. 
By the uniform density condition, there exists $K>1$ such that for some $C>0$
\[
1\leq\#\{\lambda : D_\lambda^\varphi\cap D(z,1-1/K) \neq \emptyset\}\leq C\qquad \textrm{for all $z\in\D$}.
\]
There is no restriction in assuming that $K=2$, and equivalently, that
%Replacing if necessary the dyadic squares $Q_{n,k}$ by  corresponding 
%$K$-adic squares, we can assume that
\[
1\leq\#\{\lambda : D_\lambda^\varphi\cap Q_{n,k}\neq \emptyset\}\leq C\qquad \textrm{for all $n\in\N$ and $k=0,\dots,2^n-1$}.
\]

In order to check the conditions of Theorem~\ref{suficient}(a), and given a Blaschke sequence $Z$, 
let us see first that
\[
w_{n,k}(\Lambda(\varphi),N_{n,k})\simeq N_{n,k} \log\frac 1{\varphi(2^{-n})}\ .
\]
Let $\lambda$ be such that $D_\lambda^\varphi \cap Q_{n,k} \neq \emptyset$. Take $a_1,\dots, a_N\in D_\lambda^\varphi$; 
then $\rho(a_j,\lambda)\leq \varphi(2^{-n})$ and therefore, if there is only one such disk overlapping with $Q_{n,k}$,
\[
\sum_{j=1}^{N_{n,k}}\log \frac 1{\rho(z,a_j)}\succeq N_{n,k} \log \frac 1{\varphi(2^{-n})}\ .
\]
If there is a finite number $C$ of such disks, put $N_{n,k}/C$ points in each disk, and the same result will hold.

On the other hand
\begin{align*}
\frac 1{|D_\lambda^\varphi|}\int_{D_\lambda^\varphi}\sum_{j=1}^{N_{n,k}}\log \frac 1{\rho(z,a_j)}dm(z)&\simeq
\sum_{j=1}^{N_{n,k}}\frac 1{(2^{-n}\varphi(2^{-n}))^2}\int_{D_\lambda^\varphi} \log \frac 1{\rho(z,a_j)}dm(z)\\
&\preceq\sum_{j=1}^{N_{n,k}}\frac 1{(\varphi(2^{-n}))^2}\int_{D(0,\varphi(2^{-n}))} \log \frac 1{|z|}dm(z)\\
&=\sum_{j=1}^{N_{n,k}} \log\frac 1{\varphi(2^{-n})}=N_{n,k} \log\frac 1{\varphi(2^{-n})}\ ,
\end{align*}
which proves the reverse estimate.

Let $N_n=\sum_{k=0}^{2^n-1} N_{n,k}$, $\gamma_n\in (0,1)$ to be determined later, and
$L_n=[(1-\gamma_n) N_n]$. Let $\mathcal Q_n^c$ be the set of indices $(n,k)$ corresponding 
to the $L_n$ dyadic squares $Q_{n,k}$ with the largest values of $N_{n,k}$.  By definition
\[
\sum_{(n,k)\in \mathcal Q_n^c} N_{n,k}\geq L_n\ .
\]
Call $\mathcal Q_n$ the remaining indices $(n,k)$ and define $\mathcal  Q=\cup_n\mathcal Q_n$.
Then
\begin{align*}
\sum_{(n,k)\in\mathcal Q}  2^{-n} w_{n,k}(\Lambda(\varphi), N_{n,k}) & \simeq
\sum_n  2^{-n}  \bigl(\log\frac 1{\varphi(2^{-n})}\bigr) \bigl( \sum_{k:(n,k)\in\mathcal Q} N_{n,k}\bigr)\ .
\end{align*}
Since
\[
\sum_{k:(n,k)\in\mathcal Q} N_{n,k}=N_n- \sum_{k:(n,k)\in\mathcal Q^c} N_{n,k}\leq N_n-L_n\leq N_n \gamma_n\ ,
\]
condition \eqref{invul} is now equivalent to
\[
\sum_n \gamma_n \bigl(\log\frac 1{\varphi(2^{-n})}\bigr) 2^{-n} N_n<\infty\ .
\]

The hypothesis, as stated in Remark~\ref{sum}, implies
\[
\liminf_{n\to\infty}\ \bigl(\log\frac 1{\varphi(2^{-n})}\bigr) 2^{-n} N_n=0\ ,
\]
since otherwise $2^{-n} N_n\succeq (\log  1/\varphi(2^{-n}))^{-1}$ and the Blaschke condition would be violated.
In particular, there exists a subsequence such that
\[
\sum_j (\log\frac 1{\varphi(2^{-n_j})}\bigr) 2^{-n_j} N_{n_j}<\infty\ .
\]
Define
\[
\gamma_n=
\begin{cases}
1-1/N_n\quad &\textrm{if $n=n_j$}\\
\left(\log  1/\varphi(2^{-n})\right)^{-1} &\textrm{if $n\neq n_j$.}
\end{cases}
\]
Then \eqref{invul} holds:
\[
\sum_n\gamma_n \bigl(\log\frac 1{\varphi(2^{-n})}\bigr) 2^{-n} N_n\leq
\sum_j \bigl(\log\frac 1{\varphi(2^{-n_j})}\bigr) 2^{-n_j} N_{n_j}+
\sum_n 2^{-n} N_n<\infty\ .
\]
To prove \eqref{HLzeta} in Theorem~\ref{suficient} we use an argument as in \eqref{sumLn}.
Here
\[
\sum_{n=0}^\infty \sum_{k: (n,k)\in\mathcal Q_n^c}\left(\frac{2^{-n}}{2^{-n}+k2^{-n}}\right)^2\geq
\sum_n\frac 1{L_n}\geq \sum_n\frac 1{(1-\gamma_n) N_n}\geq\sum_j 1=\infty\ ,
\]
as desired.
\end{proof}

\section{Uniformly dense disks for subharmonic funcitons}\label{S5}

In this section we show that Theorem~\ref{ud}, with a different proof, can be extended to 
the class 
\[
\mathcal S\NN=\bigl\{u:\D\longrightarrow \mathbb R\ \ \textrm{subharmonic with}\
\sup_{r<1}\int_0^{2\pi} u^+(r e^{i\theta})\; d\theta<\infty\bigr\} \ .
\]

A set $\Lambda\subset \D$ is called a \emph{determination set} for $\mathcal S\NN$ if
$\sup_\Lambda u=\sup_{\D} u$ for all $u\in \mathcal S\NN$. 

\begin{theorem}
A uniformly dense family of disks $\Lambda(\varphi)$ is a determination set for
$\mathcal S\NN$ if and only if \eqref{LyubSeip} holds.
\end{theorem}

\begin{proof}

The necessity of \eqref{LyubSeip} is contained in Theorem~\ref{ud}, 
since $\log|f|\in\mathcal S\NN$ whenever $f\in\NN$.

Assume now that \eqref{LyubSeip} holds. 
Let
$u\in\mathcal S\NN$ be such that $\sup_{\Lambda(\varphi)}u \leq 0$. We want to prove
that $u(p)\leq 0$ for all $p\notin\Lambda(\varphi)$. 

%Denote 
%\[
%\Omega(\Lambda,\varphi)=\D\setminus\overline{\Lambda(\varphi)}\ .
%\]
Let $R_n=1-K^{-n}$, where $K>1$ will be chosen
later on, and consider the domains
\[
\Omega_n(p,\Lambda,\varphi)=D(p,R_n)\setminus  
\bigcup_{\la : D_\la^\varphi\subset D(p,R_n)} \overline{D_\la^{\varphi}}\ .
\] 
Let
$\omega(A; p,\Omega)$ denote the harmonic measure at $p$ of a set
$A\subset\partial\Omega$, and let $\phi_p$ denote the automorphism of $\D$
exchanging $p$ and $0$. The subharmonicity of $u^+$ gives then
\begin{align*}
u^+(p)&\leq \int_{\partial\Omega_n(p,\Lambda,\varphi)} u^+(\zeta)\;
d\omega(\zeta;p,\Omega_n(p,\Lambda,\varphi))
=\int_{\partial D(p,R_n)} u^+(\zeta)\;
d\omega(\zeta;p,\Omega_n(p,\Lambda,\varphi))\\
&=\int_{|\zeta|=R_n} (u^+\circ\phi_p)(\zeta)\;
d\omega(\zeta;0,\phi_p(\Omega_n(p,\Lambda,\varphi)))\ .
\end{align*}

First observe that the harmonic measure in $\phi_p(\Omega_n(p,\Lambda,\varphi)))$ 
can be estimated by the harmonic measure of a domain 
$\Omega_n(0,\tilde\Lambda,\psi)$, where $\tilde\Lambda$ is uniformly dense and $\psi$ is a
non-decreasing, continuous function bounded by some constant less than 1 satisfying \eqref{LyubSeip}.
To see this let $\tilde\Lambda=\phi_p(\Lambda)$, consider the hyperbolic rings
\[
A_n=\{z\in\D : R_{n-1}\leq\rho(z,p)< R_n\}
\]
and take $\psi$ non-decreasing, continuous, and
such that $\psi(R_n)=\min_{A_n}\varphi$. Then
$D_{\phi_p(\lambda)}^\psi\subset D_{\phi_p(\lambda)}^\varphi$
and therefore $\omega(A;0,\phi_p(\Omega_n(p,\Lambda,\varphi)))\leq
\omega(A;0,\Omega_n(0,\tilde\Lambda,\psi))$ 
for any $A\subset\{|z|=R_n\}$.
Notice also that $\min_{A_n}\varphi$ is
attained for $z$ with
\[
1-|z|=\frac{(1-R_n)(1-|p|)}{1+R_n|p|}\leq
(1-R_n)(1-|p|)= 
K^{-n}(1-|p|) \ ,
\]
and therefore
\[
\int_0^1\frac{dt}{t\log(1/\psi(t))}\simeq
\sum_{n=1}^\infty \frac 1{\log(1/\psi(K^{-n}))}\leq
\sum_{n=1}^\infty \frac 1{\log(1/\varphi(K^{-n}(1-|p|))}\simeq
\int_0^1\frac{dt}{t\log(1/\varphi(t))}\ .
\]

We have thus
\[
u^+(p)\leq \int_{|\zeta|=R_n} (u^+\circ\phi_p)(\zeta)\;
d\omega(\zeta;0,\Omega_n(0,\tilde\Lambda,\psi))\ .
\]
As mentioned in Remark~\ref{sum}, the 
hypothesis implies
$\omega(\partial\D;0,\D\setminus\overline{\tilde\Lambda(\psi)})=0$. 
In order to see that the previous integrals tend to zero 
we need a slight refinement of Theorem 1 in \cite{OrSe}.
Let $d\sigma_n=d\theta/(2\pi R_n)$ denote the normalized Lebesgue measure
in $|z|=R_n$.

\begin{lemma}\label{mesuraharmonica} Given a uniformly dense sequence
$\Lambda$ and a non-decreasing continuous function $\varphi$ satisfying \eqref{LyubSeip},
there exist $R_n<1$ with $\lim_n R_n=1$, and 
$\epsilon_n>0$ with $\lim_n \epsilon_n=0$ such that 
$\omega(I;0,\Omega_n(0,\Lambda,\varphi))\leq \epsilon_n \sigma_n(I)$
for all intervals $I\subset\{\zeta : |\zeta|=R_n\}$.
\end{lemma}

Once this lemma is proved, the above estimate yields
\[
u^+(p)\leq \epsilon_n \int_{|\zeta|=R_n} 
(u^+\circ\phi_p)(\zeta)\; d\sigma_n(\zeta) \leq \frac{\epsilon_n}{R_n} 
\sup_{r<1}\int_0^{2\pi} (u^+\circ\phi_p)(r e^{i\theta})\; \frac{d\theta}{2\pi}\ ,
\]
and letting $n\to\infty$ we obtain $u^+(p)\leq 0$, as desired.
\end{proof}

\begin{proof}[Proof of Lemma~\ref{mesuraharmonica}]
We prove this by induction. We drop the superindex in $D_\la^{\varphi}$ and
denote $\Omega_n(0,\Lambda,\varphi)$ simply by $\Omega_n$.

There is no restriction in assuming that there
are no $D_\la$ in $D(0,R_1)$; thus $\omega(I;0,\Omega_1))\leq  |I|$ for all 
intervals $I\subset\{\zeta : |\zeta|=R_1\}$.

We have
\[
\omega(I;0,\Omega_n)=\int_{|z|=R_{n-1}} P(z\to I)\;
 d\omega(z;0,\Omega_{n-1})\ ,
\]
where $P(z\to I)$ denotes the probability that a Brownian motion starting at $z$
exits $\Omega_n$ through $I$. The hypothesis of induction gives then
\begin{equation}\label{HI}
\omega(I;0,\Omega_n)\leq\epsilon_{n-1} \int_{|z|=R_{n-1}} P(z\to I)\; d\sigma_{n-1}(z)\ .
\end{equation}
In the estimate of $P(z\to I)$ we use the uniform density of $\Lambda$: there exist
$\delta\in (0,1)$ and $K>1$ (independent of $n$) such that for each $z\in\{|z|=R_{n-1}\}$
there is $\la\in \Lambda$ with $D_\la\in D(0,R_n)\setminus D(0,R_{n-1})$ and
$\rho(z,\la)\leq \delta$. Then
\[
P(z\to I)\leq \omega(I;z, D(0,R_n)\setminus D_{\lambda})\ .
\]
This harmonic measure can be estimaded by comparing with an explicit harmonic function.
Let $\Psi_n(z)=z/R_n$, which sends $D(0,R_n)$ to $\D$, and let
$P_z^{(n)}(\zeta)$ denote the Poisson kernel in $D(0,R_n)$.
Let $\la^{(n)}\in\D$, $r_\la^{(n)}\in (0,1)$
be such that $\Psi_n(D_\lambda)=\Psi_n(D(\lambda,r_\lambda))=D( \la^{(n)},r_\la^{(n)})$, and define
the harmonic function on $D(0,R_n)\setminus D_\lambda$
\[
F_\la(z,I)=\omega(I;z,D(0,R_n))-
%\int_I P_z^n(\zeta)\; d\sigma_n(\zeta)-
\left(\inf_{w\in\partial D_\lambda} \int_I P_w^{(n)}(\zeta) d\sigma_n(\zeta)
\right)
\frac{\log\rho(\la^{(n)},\Psi_n(z))}{\log r_\lambda^{(n)}}\ .
\]
It is clear that $\omega(I;\eta, D(0,R_n)\setminus D_\lambda)\leq F_\la(\eta,I)$ for $\eta$
in the boundary of $D(0,R_n)\setminus D_\lambda$, and therefore in all $D(0,R_n)\setminus D_\lambda$. Hence
\[
P(z\to I)\leq F_\la(z,I)\ .
\]

We want to give an estimate of $F_\la(z,I)$ independent of $\lambda$. Using that
$\omega(I;z,D(0,R_n))=\int_I P_z^{(n)}(\zeta)\; d\sigma_n(\zeta)$,  we can write
\[
F_\la(z,I)=\left(\int_I P_z^{(n)}(\zeta)\; d\sigma_n(\zeta)\right)
\left(1-\frac{\inf_{w\in\partial D_\lambda} \int_I P_w^{(n)}(\zeta)
d\sigma_n(\zeta)}{\int_I P_z^{(n)}(\zeta)\; d\sigma_n(\zeta)} 
\frac{\log\rho(\la^{(n)},\Psi_n(z))}{\log r_\lambda^{(n)}}\right)\ .
\]
Since $\lim_{|\la|\to 1} r_\la=0$
and $\rho(z,\la)\leq \delta$, 
by the Harnack's estimates (or by a direct computation), there exists $c>0$ such that
\[
\frac{\inf_{w\in\partial D_\lambda} \int_I P_w^{(n)}(\zeta) d\sigma_n(\zeta)}{\int_I P_z^{(n)}(\zeta)
\; d\sigma_n(\zeta)}\geq c\ .
\]
Also, there exists $\delta'>0$ such that for $n$ big enough
$ \rho(\la^{(n)},\Psi_n(z))\leq \delta'$.
With this and the fact that $K^{-n-1}\leq 1-|\la| < K^{-n}$ we deduce that there exists some 
$C>0$ such that
\[
P(z\to I)\leq \left(\int_I P_z^{(n)}(\zeta)\; d\sigma_n(\zeta)\right)
\left(1-\frac{C}{\log 1/\varphi(K^{-n})}\right)\ .
\]
From \eqref{HI} we have therefore
\begin{align*}
\omega(I;0,\Omega_n)&\leq\epsilon_{n-1} \left(1-\frac{C}{\log 1/\varphi(K^{-n})}\right)
\int_I \int_{|z|=R_{n-1}} P_z^{(n)}(\zeta)\; d\sigma_{n-1}(z) \; d\sigma_n(\zeta)\\
&=\epsilon_{n-1} \left(1-\frac{C}{\log 1/\varphi(K^{-n})}\right) \sigma_n(I)\ .
\end{align*}
Defining 
\[
\epsilon_n=\epsilon_{n-1} \left(1-\frac{C}{\log 1/\varphi(K^{-n})}\right)
=\prod_{j=1}^n \left(1-\frac{C}{\log 1/\varphi(K^{-j})}\right)\ 
\]
and using Remark~\ref{sum} we obtain the stated properties.
\end{proof}

\section{Appendix.}

Here we give a direct proof that a determination set for $\NN$ is
a Hayman-Lyons set. According to \cite[Corollary 2]{Ga} this implies
$|NT(\Lambda)|=2\pi$, and therefore $\Lambda$ is determination set for
$\H^\infty$.

Let $\delta\in (0,1)$ and let $\Lambda_0\subset\Lambda$ be maximal among the
subsequences of $\Lambda$ such that $\rho(\lambda,\lambda')\geq \delta$ for all
$\lambda, \lambda'\in\Lambda_0$, $\lambda\neq\lambda'$. We want to prove that
\[
\sum_{\lambda\in \Lambda_0}(1-|\lambda|) P_\lambda(\zeta)=
\infty\qquad\textrm{for all $\zeta\in\T$}.
\]
There is no loss of generality in reducing ourselves to the case $\zeta=1$.
Also, we restrict our attention to those $\lambda\in\Lambda_0$ with
$P_\lambda(1)\geq 1$, and denote by $\tilde\Lambda_0$ the subsequence made with
such points. Notice that $\Lambda\cap\{z: P_z(1)\leq 1\}$ cannot be a
determination sequence for $\NN$ anyway, as the function $f\in\NN$ with
$\log|f(z)|=P_z(1)$ shows.  Thus,  let us assume that
\[
\sum_{\lambda\in \tilde\Lambda_0}(1-|\lambda|) P_\lambda(1)<\infty\ 
\]
and see that there exists
$f\in\NN\setminus\H^\infty$ with $\sup_\Lambda|f|<\infty$.

Consider the sequence $Z$ consisting of the points
$\lambda\in\tilde\Lambda_0$, with multiplicity $[P_\lambda(1)]$. By assumption
$Z$ is a Blaschke sequence, and therefore, for any $C>0$, there exists
$f\in\NN$ such that  
\[  
\log|f(z)|=\log|B_Z(z)|+C\; P_z(1)\ .  
\]  
As seen in Example 1, such $f$ cannot be bounded in the whole disk. 

We want to choose $C$ so that $\sup_\Lambda|f|<\infty$. By construction, we only
need to consider $\lambda\notin\tilde\Lambda_0$. We separate two cases:
\begin{itemize}
\item[(i)] If $P_\lambda(1)\leq 2$ obviously $\log|f(\lambda)|\leq 2C$. 
\item[(ii)] If $P_\lambda(1)> 2$ there exists $\lambda_0\in \tilde\Lambda_0$ such
that $\rho(\lambda,\lambda_0)\leq\delta$. By Harnack's inequalities we
obtain:
\[
\log |f(\lambda)|\leq\log\rho(\lambda,\lambda_0)^{[P_{\lambda_0}(1)]}+C \;
P_\lambda(1)\leq
\frac 12(\log\delta) P_{\lambda_0}(1)+C \bigl(\frac{1+\delta}{1-\delta}\bigr)
P_{\lambda_0}(1)\ .
\]
Choosing $C=\frac 12\frac{1-\delta}{1+\delta}\log\frac 1\delta$ we see that
in this case $\log |f(\lambda)|\leq 0$, as desired.

\end{itemize}

\end{document}